\documentclass[11pt]{article}
\usepackage{a4wide,amssymb,euscript,stmaryrd,epsfig}

\usepackage{amssymb,amsmath,amsthm}

\usepackage{epsfig}
\usepackage{amsfonts}
\usepackage{latexsym}
\usepackage{xspace}
\usepackage{graphicx}
\usepackage[usenames]{color}
\usepackage{subfigure}
\usepackage{url}

\def\nnew{}

\newtheorem{theorem}{Theorem}[section]
\newtheorem{lemma}[theorem]{Lemma}
\newtheorem{proposition}[theorem]{Proposition}

\newtheorem{definition}[theorem]{Definition}

\newtheorem{remark}[theorem]{Remark}

\newcommand{\supp}{\operatorname{supp}}

\newcommand{\R}{\mathbb R}

\newcommand{\N}{\mathbb N}

\title{A Convergent Overlapping Domain Decomposition Method for Total Variation Minimization}

\author{Massimo Fornasier\thanks{Johann Radon Institute for Computational and Applied Mathematics (RICAM),
Austrian Academy of Sciences, Altenbergerstrasse 69, A-4040, Linz, Austria Email: {\tt massimo.fornasier@oeaw.ac.at}} \and Andreas Langer\thanks{Johann Radon Institute for Computational and Applied Mathematics (RICAM),
Austrian Academy of Sciences, Altenbergerstrasse 69, A-4040, Linz, Austria Email: {\tt andreas.langer@oeaw.ac.at}} \and Carola-Bibiane Sch\"onlieb\thanks{Department of Applied Mathematics and Theoretical Physics (DAMTP),
Centre for Mathematical Sciences,
Wilberforce Road,
Cambridge CB3 0WA,
United Kingdom.Email: {\tt c.b.s.schonlieb@damtp.cam.ac.uk} } }

\begin{document}

\graphicspath{{./graphics/}}

\maketitle

\pagestyle{myheadings}
\thispagestyle{plain}
\markboth{M. FORNASIER, A. LANGER, AND C.-B. SCH\"ONLIEB}{OVERLAPPING DOMAIN DECOMPOSITION METHODS FOR TV-MINIMIZATION}

\begin{abstract}
This paper is concerned with the analysis of convergent sequential and parallel overlapping domain decomposition methods for the minimization of functionals formed by a discrepancy term with respect to data and {\nnew a total variation constraint}. To our knowledge, this is the first successful attempt of addressing such strategy for the nonlinear, nonadditive, and nonsmooth problem of total variation minimization. We provide several numerical experiments, showing the successful application of the algorithm for the restoration of 1D signals and 2D images in interpolation/inpainting problems respectively, and in a compressed sensing problem, for recovering piecewise constant medical-type images from partial Fourier ensembles.
\end{abstract}

\noindent {\bf Key words}: 
Domain decomposition method, nonsmooth convex optimization, parallel computation, discontinuous solutions, total variation minimization, $\ell_1$-minimization, image and signal processing\\

\noindent {\bf AMS subject classifications.} 
65K10, 
65N55  
65N21, 
65Y05  
90C25, 
52A41, 
49M30, 
49M27, 
68U10  

\section{Introduction}

In concrete applications, e.g., for image processing, one might be interested to recover at best a digital image provided only partial linear or nonlinear measurements, possibly corrupted by noise.
Given the observation that natural and man-made images can be characterized by a relatively small number of edges and extensive relatively uniform parts, one may want to help the reconstruction by imposing that the interesting solution is the one which matches the given data and has also a few discontinuities localized on sets of lower dimension.

In the context of {\it compressed sensing} \cite{cand,carotaXX,cataXX,do04}, it has been clarified {\nnew that the minimization of $\ell_1$-norms occupies a fundamental role for the promotion of sparse solutions}. This understanding furnishes an important interpretation of {\it total variation minimization}, i.e., the minimization of the $L^1$-norm of derivatives \cite{ROF}, as a regularization technique for image restoration. The problem can be modelled as follows; let $\Omega \subset \mathbb R^d$, for $d=1,2$ be a bounded open set  with Lipschitz boundary, and $\mathcal H = L^2 (\Omega)$. For $u \in L_{loc}^1(\Omega)$
$$
V(u,\Omega) := \sup \left  \{ \int_\Omega u \operatorname{div} \varphi~dx: \varphi \in \left [ C^1_c(\Omega) \right ]^d, \| \varphi \|_\infty \leq 1 \right \}
$$
is the variation of $u${\nnew . Further, }$u \in BV(\Omega)$, the space of bounded variation functions \cite{AFP,EvGa},  if and only if $V(u,\Omega) < \infty$. In this case, we denote $|D(u)|(\Omega) =V(u,\Omega)$. If $u \in W^{1,1}(\Omega)$ (the Sobolev space of $L^1$-functions with $L^1$-distributional derivatives), then $|D(u)|(\Omega) = \int_{\Omega} | \nabla u | \, d\,x$.
We consider as in \cite{ChL,Ve01} the minimization in $BV(\Omega)$ of the functional 
\begin{equation}\label{functotal1}
\mathcal{J}(u):= \| T u - g \|_2^2 + 2 \alpha \left|D(u)\right|(\Omega),
\end{equation} 
where $T:L^2(\Omega) \to L^2(\Omega)$ is a bounded linear operator, $g \in L^2(\Omega)$ is a datum, and $\alpha>0$ is a fixed {\it regularization parameter} \cite{enhane96}.
Several numerical strategies to perform efficiently total variation minimization have been proposed in the literature. {\nnew Without claiming of being exhaustive}, we list a few of the relevant methods, ordered by their chronological appearance:

(i) the linearization approach of Vogel et al. \cite{dovo} and of Chambolle and Lions \cite{ChL} by iteratively re-weighted least squares, see also \cite{DDFG} for generalizations and refinements in the context of compressed sensing;

(ii) the primal-dual approach of Chan et al. \cite{chgomu};

(iii) variational approximation via locally quadratic functionals as in the work of Vese et al. \cite{AK02,Ve01};

(iv) iterative thresholding algorithms based on projections onto convex sets as in the work of Chambolle \cite{Ch} as well as in the work of Combettes and Wajs \cite{CW} and Daubechies et al. \cite{dateve06};

(v)  iterative minimization of the Bregman distance as in the work of Osher et al. \cite{breg} {\nnew (also notice the very recent Bregman split approach \cite{GoOs})};

{\nnew (vi) graph cuts \cite{ChDa,DaSig} for the minimization of \eqref{functotal1} with $T=Id$ and an anisotropic total variation;} 

(vii)  the approach proposed by Nesterov \cite{ne05} and its modifications by Weiss et al. \cite{WBA}.

These approaches differ significantly, and they provide a convincing view of the interest this problem has been able to generate and of his applicative impact. 
However, because of their iterative-sequential formulation, none of the mentioned methods is able to address in real-time, or at least in an acceptable computational time, extremely large problems, such as 4D imaging (spatial plus temporal dimensions) from functional magnetic-resonance in nuclear medical imaging, astronomical imaging or global terrestrial seismic tomography.
For such large scale simulations we need to address methods which allow us to reduce the problem to a finite sequence of sub-problems of a more manageable size, perhaps {\nnew computable} by one of the methods listed above. With this aim we introduced subspace correction and domain decomposition methods  both for $\ell_1$-norm and total variation minimizations \cite{fo07,FS,sc09}. We address the interested reader to the broad literature included in \cite{FS} for an introduction to domain decompositions methods both for PDEs and convex minimization.

\subsection{Difficulty of the problem} Due to the nonsmoothness and nonadditivity of the total variation with respect to a nonoverlapping domain decomposition (note that the total variation of a function on the whole domain equals the sum of the total variations on the subdomains plus the size of the jumps at the interfaces \cite[formula (3.4)]{FS}), one encounters additional difficulties in showing convergence of such decomposition strategies to global minimizers. In particular, we stress very clearly that well-known approaches as in {\nnew \cite{CC,ChMa,TT,TX}} are not directly applicable to this problem, because either they do address additive problems or smooth convex minimizations, which is {\it not} the case of total variation minimization. 
 Moreover the interesting  solutions may be discontinuous, e.g., along curves in 2D. These discontinuities may cross the interfaces of the domain decomposition patches. Hence, the crucial difficulty is the correct numerical treatment of interfaces, with the preservation of crossing discontinuities and the correct matching where the solution is continuous instead, see \cite[Section 7.1.1]{FS}.

The work \cite{FS} was particularly addressed to {\it nonoverlapping} domain decompositions {\nnew $\Omega_1 \cup \Omega_2 \subset \Omega \subset \bar \Omega_1 \cup \bar \Omega_2$} and $\Omega_1 \cap \, \Omega_2 = \emptyset$. Associated to the decomposition define $V_i=\{ u \in L^2(\Omega): \operatorname{supp}(u) \subset \Omega_i \}$, for $i=1,2$; note that $L^2(\Omega) = V_1 \oplus V_2$. With this splitting we wanted to minimize $\mathcal J$ by suitable instances of the following alternating algorithm:  Pick an initial $V_1 \oplus V_2 \ni  u_1^{(0)}+ u_2^{(0)} : = u^{(0)}$, for example $u^{(0)}=0$, and iterate
$$
\left \{ 
\begin{array}{ll}
u_1^{(n+1)} \approx \arg \min_{v_1 \in V_1}  \mathcal  J(v_1 +u_2^{(n)}) &\\
u_2^{(n+1)} \approx  \arg \min_{v_2 \in V_2} \mathcal J(u_1^{(n+1)} + v_2) &\\
u^{(n+1)}:=u_1^{(n+1)} + u_2^{(n+1)}.
\end{array}
\right.
$$ 
In \cite{FS} we proposed an implementation of this algorithm which is guaranteed to converge and {\nnew to decrease the objective energy $\mathcal J$ monotonically.} We could prove its convergence to minimizers of $\mathcal J$ only under technical conditions on the interfaces of the subdomains. However, in our numerical experiments, the algorithm seems always converging robustly to the expected minimizer. This discrepancy between theoretical analysis and numerical evidences motivated our investigation on {\it overlapping} domain decompositions. The hope was that the redundancy given by overlapping patches and the avoidance of boundary interfaces could allow {\nnew for a technically easier theoretical analysis.}

\subsection{Our approach, results, and technical issues} In this paper we show how to adapt our previous algorithm \cite{FS} to the case of an {\it overlapping} domain decomposition. {\nnew The setting of an overlapping domain decomposition eventually provides us with a framework in which }we successfully prove its convergence to minimizers of $\mathcal J$, both in its sequential and parallel forms. Let us stress that to our knowledge this is the first method which addresses a domain decomposition strategy for total variation minimization with a formal theoretical justification {\nnew of convergence}.  It is important to mention that there are other very recent attempts of addressing domain decomposition methods for total variation minimization with successful numerical results \cite{JM}.

Our analysis is performed {\nnew for a discrete approximation of the continuous functional \eqref{functotal1}, for ease again denoted $\mathcal J$ in \eqref{functotalfin}.} Essentially we approximate functions $u$ by their sampling on a regular grid and their gradient $D u$ by finite differences $\nabla u$. It is well-known that such a discrete approximation $\Gamma$-converges  to the continuous functional (see \cite{br02}). In particular, discrete minimizers of \eqref{functotalfin}, interpolated by piecewise linear functions, converge in weak-$*$-topology of $BV$ to minimizers of the functional \eqref{functotal1} in the continuous setting. 
Of course, when dealing with numerical solutions, only the discrete approach matters together with its approximation properties to the continuous problem. However, the need of {\nnew working in the discrete setting is not only} practical, it is also topological. In fact bounded sets in $BV$ are (only) weakly-$*$-compact, and this property is fundamental for showing that certain sequences have converging subsequences. Unfortunately, the weak-$*$-topology of $BV$ is ``too weak'' {\nnew for our purpose of proving convergence of the domain decomposition algorithm}; for instance, the trace on boundary sets is {\it not} a continuous operator with respect to {\nnew this} topology. This difficulty can be avoided, for instance, by $\Gamma$-approximating the functional \eqref{functotal1} by means of quadratic functionals (as in \cite{AK02,ChL,Ve01}) and working with the topology of $W^{1,2}(\Omega)$, the Sobolev space of $L^2$-functions with $L^2$-distributional {\nnew first} derivatives. However, this strategy changes the singular nature of the problem which makes it both interesting and difficult. Hence, the discrete approach has the virtues of being practical for numerical implementations, {\nnew of correctly} approximating the continuous setting, and {\nnew of} retaining the major features which makes the problem interesting. Note further that in the discrete setting where topological issues are not a concern anymore, also the dimension $d$ can be arbitrary, contrary to the continuous setting where the dimension $d$ has to be linked to boundedness properties of the operator $T$, see \cite[property H2, pag. 134]{Ve01}.
For ease of presentation, and in order {\nnew to avoid} unnecessary technicalities, we limit our analysis to splitting the problem into two subdomains $\Omega_1$ and $\Omega_2$. This is by no means a restriction. The generalization to multiple domains comes quite natural in our specific setting, see also \cite[Remark 5.3]{FS}. When dealing with discrete subdomains $\Omega_i$, for technical reasons, we will require a certain splitting property for the total variation, i.e.,
\begin{equation}
\label{splitdom}
|\nabla u|(\Omega) = |\nabla u|_{\Omega_1}|(\Omega_1) + c_1(u|_{(\Omega_2 \setminus \Omega_1) \cup \Gamma_1}), \quad |\nabla u|(\Omega) = |\nabla u|_{\Omega_2}|(\Omega_2) + c_2(u|_{(\Omega_1\setminus \Omega_2) \cup \Gamma_2}),
\end{equation}
where $c_1$ and $c_2$ are suitable functions which depend only on the restrictions $u|_{(\Omega_2 \setminus \Omega_1) \cup \Gamma_1}$ and $u|_{(\Omega_1 \setminus \Omega_2) \cup \Gamma_2}$ respectively, see \eqref{tvdecomp} (symbols and notations are clarified once for all in the following section). Note that this formula is the discrete analogous of  \cite[formula (3.4)]{FS} in the continuous setting.
{\nnew The simplest examples} of discrete domains with such a property are {\nnew discrete $d$-dimensional rectangles ({\it $d$-orthotopes})}. Hence, for ease of presentation, we will assume to work with $d$-orthotope domains, also noting that such decompositions are already sufficient for any practical use in image processing, and stressing that the results can be generalized also to subdomains with different shapes as long as \eqref{splitdom} is satisfied.

\subsection{Organization of the work}
The paper is organized as follows. In Section \ref{notations} we collect the relevant notations and symbols for the paper. Section \ref{Problem} introduces the problem and the overlapping domain decomposition algorithm which we want to analyze. In Section \ref{localproblem} we address the approximate solution of the local problems defined on the subdomains $\Omega_i$ and we show how we interface them, by means of a suitable Lagrange multiplier. Section \ref{Convsec} and Section \ref{Convsec2} are concerned with the convergence of the sequential and parallel forms of the algorithm. In particular, in Section  \ref{Convsec} we provide a characterization of minimizers by a discrete representation of the subdifferential of $\mathcal J$. This characterization is used in the convergence proofs in order to check the reached optimality. The final Section \ref{Numersec} provides a collection of applications and numerical examples.

\section{Notations}\label{notations}
Let us fix the main notations. Since we are interested in a discrete setting we define the {\it discrete $d$-orthotope} $\Omega = \{x_1^1 < \ldots < x_{N_1}^1\} \times \ldots \times \{x_1^d < \ldots < x_{N_d}^d\}\subset\mathbb{R}^d$, $d\in\mathbb{N}$ and the considered function spaces are $\mathcal H = \R^{N_1\times N_2 \times \ldots \times N_d}$, where $N_i\in\N$ for $i=1,\ldots,d$. For $u\in\mathcal{H}$ we write $u=u(x_{\mathbf{i}})_{\mathbf{i}\in\mathcal{I}}$ with
$$
\mathcal{I}:=\prod_{k=1}^d\{1,\ldots,N_k\}
$$
and
$$
u(x_{\mathbf{i}}) = u(x_{i_1}^1, \ldots , x_{i_d}^d)
$$
where $i_k\in\{1,\ldots,N_k\}$ and $(x_{\mathbf{i}})_{\mathbf{i}\in\mathcal{I}} \in \Omega$. Then we endow $\mathcal{H}$ with the norm
$$
\|u\|_{\mathcal{H}}=\|u\|_2=\left(\sum_{\mathbf{i}\in\mathcal{I}}|u(x_{\mathbf{i}})|^2\right)^{1/2}=\left(\sum_{x\in\Omega}|u(x)|^2\right)^{1/2}.
$$
We define the scalar product of $u,v \in \mathcal{H}$ as
$$
\langle u,v \rangle_{\mathcal{H}}= \sum_{\mathbf{i}\in\mathcal{I}}u(x_{\mathbf{i}})v(x_{\mathbf{i}})
$$
and the scalar product of $p,q \in \mathcal{H}^d$ as
$$
\langle p,q \rangle_{\mathcal{H}^d}= \sum_{\mathbf{i}\in\mathcal{I}}\langle p(x_{\mathbf{i}}),q(x_{\mathbf{i}}) \rangle_{\R^d} 
$$
with $\langle y,z \rangle_{\R^d} = \sum_{j=1}^{d} y_j z_j$ for every $y=(y_1,\ldots,y_d)\in\R^d$ and $z=(z_1,\ldots,z_d)\in\R^d$. 
We will consider also other norms, in particular 
$$
\|u\|_p=\left(\sum_{\mathbf{i}\in\mathcal{I}}|u(x_{\mathbf{i}})|^p\right)^{1/p}, \quad 1 \leq p < \infty,
$$
and
$$
\| u \|_\infty = \sup_{\mathbf{i}\in\mathcal{I}}|u(x_{\mathbf{i}})|.
$$
We denote the discrete gradient $\nabla u$ by
$$
(\nabla u)(x_{\mathbf{i}}) = ((\nabla u)^1(x_{\mathbf{i}}),\ldots,(\nabla u)^d(x_{\mathbf{i}}))
$$
with 
$$
(\nabla u)^j(x_{\mathbf{i}})=
\begin{cases}
u(x^1_{i_1},\ldots,x^j_{i_j+1},\ldots,x^d_{i_d}) - u(x^1_{i_1},\ldots,x^j_{i_j},\ldots,x^d_{i_d}) & \text{if } i_j<N_j \\
0 & \text{if } i_j=N_j
\end{cases}
$$
for all $j=1,\ldots,d$ and for all $\mathbf{i}=(i_1,\ldots,i_d)\in\mathcal{I}$.\\
{\nnew Let $\varphi:\R \to \R$, we define for $\omega\in\mathcal{H}^d$
$$
\varphi(|\omega|)(\Omega)=\sum_{\mathbf{i}\in\mathcal{I}} \varphi(|\omega(x_{\mathbf{i}})|) = \sum_{x\in\Omega} \varphi(|\omega(x)|),
$$
where $|y|=\sqrt{y_1^2+\ldots+y_d^2}$.
In particular we define the \emph{total variation} of $u$ by setting $\varphi(s)=s$ and $\omega=\nabla u$, i.e., 
$$
|\nabla u|(\Omega):=\sum_{\mathbf{i}\in\mathcal{I}}|\nabla u (x_{\mathbf{i}})|= \sum_{x\in\Omega}|\nabla u(x)|.
$$}
%
For an operator $T$ we denote $T^*$ its adjoint. Further we introduce the \emph{discrete divergence} $\operatorname{div}: \mathcal{H}^d \to \mathcal{H}$ defined, in analogy with the continuous setting, by $\operatorname{div} = -\nabla^*$ ($\nabla^*$ is the adjoint of the gradient $\nabla$). The discrete divergence operator is explicitly given by 
\begin{eqnarray*}
(\operatorname{div} p)(x_{\mathbf{i}}) &=& \begin{cases}
p^1(x^1_{i_1},\ldots,x^d_{i_d})-p^1(x^1_{i_1-1},\ldots,x^d_{i_d}) & \textrm{ if } 1<i_1<N_1\\
p^1(x^1_{i_1},\ldots,x^d_{i_d}) & \textrm{ if } i_1=1\\
-p^1(x^1_{i_1-1},\ldots,x^d_{i_d}) & \textrm{ if } i_1=N_1
\end{cases}
\\
& & + \ldots + 
\begin{cases}
p^d(x^1_{i_1},\ldots,x^d_{i_d})-p^d(x^1_{i_1},\ldots,x^d_{i_d-1}) & \textrm{ if } 1<i_d<N_d\\
p^d(x^1_{i_1},\ldots,x^d_{i_d}) & \textrm{ if } i_d=1\\
-p^d(x^1_{i_1},\ldots,x^d_{i_d-1}) & \textrm{ if } i_d=N_d,
\end{cases}
\end{eqnarray*}
for every $p=(p^1,\ldots,p^d)\in \mathcal{H}^d$ and for all $\mathbf{i}=(i_1,\ldots,i_d)\in\mathcal{I}$.
(Note that if we considered discrete domains $\Omega$ which are not discrete $d$-orthotopes, then the definitions of gradient and divergence operators should be adjusted accordingly.) 
With these notations, we define the closed convex set 
{\nnew \begin{equation*}
K:=\left\{\operatorname{div} p: p\in \mathcal H^d, \left|p(x)\right|_\infty\leq 1\;\mbox{ for all } x \in \Omega \right\},
\end{equation*}
where $\left|p(x)\right|_\infty = \max\left \{ |p^1(x)|,\ldots,|p^d(x)|\right \}$,} and denote $P_K(u)=\operatorname{argmin}_{v\in K} \|u-v\|_2$ the \emph{orthogonal projection onto $K$}.
We will often use the symbol $1$ to indicate the constant vector with entry values $1$ and $1_{D}$ to indicate the characteristic function of the domain $D \subset \Omega$.


\section{The Overlapping Domain Decomposition Algorithm}\label{Problem}
We are interested in the minimization of the functional
\begin{equation}\label{functotalfin}
\mathcal{J}(u):= \| T u - g \|_2^2 + 2 \alpha \left|\nabla(u)\right|(\Omega),
\end{equation} 
where $T \in \mathcal L(\mathcal{H})$ is a  linear operator, $g \in \mathcal{H}$ is a datum, and $\alpha>0$ is a fixed constant. 
In order to guarantee the existence of minimizers for \eqref{functotalfin} we assume that:
\begin{itemize}
\item[(C)] $\mathcal J$ is coercive in $\mathcal{H}$, i.e., there exists a constant $C>0$ such that $\{ \mathcal J \leq C\} := \{ u \in \mathcal{H} :  \mathcal J(u) \leq C\}$ is bounded in $\mathcal{H}$.
\end{itemize}
It is well known that if $1\notin\ker(T)$ then condition (C) is satisfied, see \cite[Proposition 3.1]{Ve01}.\\

Now, instead of minimizing \eqref{functotalfin} on the whole domain we decompose $\Omega$ into two overlapping subdomains $\Omega_1$ and $\Omega_2$ such that $\Omega=\Omega_1\cup\Omega_2$, $\Omega_1\cap\Omega_2\not=\emptyset$, and \eqref{splitdom} is fulfilled. For consistency of the definitions of gradient and divergence, we assume that also the subdomains $\Omega_i$ are discrete $d$-orthotopes as well as $\Omega$, stressing that this is by no means a restriction, but only {\nnew for} ease of presentation. 
Due to this domain decomposition $\mathcal{H}$ is split into two closed subspaces $V_j=\{ u \in \mathcal H: \operatorname{supp}(u) \subset \Omega_j \}$, for $j=1,2$. Note that $\mathcal H = V_1 + V_2$ is not a direct sum of  $V_1$ and $V_2$, but just a linear sum of subspaces.
{\nnew Thus any $u \in \mathcal H$ has a nonunique representation
\begin{equation}\label{u}
u(x)=
\begin{cases}
u_1(x) & x\in \Omega_1 \setminus \Omega_2\\
u_1(x) + u_2(x) & x\in \Omega_1 \cap \Omega_2\\
u_2(x) & x\in \Omega_2 \setminus \Omega_1 
\end{cases}, \quad u_i \in V_i, \quad i=1,2.
\end{equation}}
We denote by $\Gamma_1$ the interface between $\Omega_1$ and $\Omega_2\setminus\Omega_1$ and by $\Gamma_2$ the interface between $\Omega_2$ and $\Omega_1\setminus\Omega_2$ (the interfaces are naturally defined in the discrete setting).






We introduce the trace operator of {\nnew the} restriction to a boundary $\Gamma_i$
$$
\operatorname{Tr}\mid_{\Gamma_i}: V_i \rightarrow \mathbb R^{\Gamma_i}, \ \ i=1,2
$$ 
with $\operatorname{Tr}\mid_{\Gamma_i} v_i = v_i\mid_{\Gamma_i}$, the restriction of $v_i$ on $\Gamma_i$. 
{\nnew Note that $\mathbb R^{\Gamma_i}$ is as usual the set of maps from $\Gamma_i$ to $\mathbb R$.}
The trace operator is clearly a linear and continuous operator. We additionally fix a \emph{bounded uniform partition of unity} (BUPU) $\{\chi_1, \chi_2 \} \subset \mathcal{H}$ such that
\begin{itemize}
\item[(a)] $\operatorname{Tr}\mid_{\Gamma_i}\chi_i = 0$ for $i=1,2$,
\item[(b)] $\chi_1 + \chi_2 = 1$,
\item[(c)] $\supp \chi_i \subset \Omega_i$ for $i=1,2$,
\item[(d)] \nnew{$\max\{||\chi_1||_\infty,||\chi_2||_\infty\}=\kappa<\infty.$}
\end{itemize} 
We would like to solve
$$
\operatorname{argmin}_{u\in\mathcal{H}}\mathcal{J}(u)
$$  
by picking an initial $V_1+ V_2 \ni \tilde{u}_1^{(0)} + \tilde{u}_2^{(0)}:=u^{(0)}\in\mathcal{H}$, e.g., $\tilde u_i^{(0)}=0, i=1,2$, and iterate
\begin{equation}\label{schw_sp}
\left \{ 
\begin{array}{ll}
u_1^{(n+1)} \approx \operatorname{argmin}_{\stackrel{v_1 \in V_1}{\operatorname{Tr}\mid_{\Gamma_1}v_1=0}}  \mathcal  J(v_1 +\tilde{u}_2^{(n)}) &\\
u_2^{(n+1)} \approx  \operatorname{argmin}_{\stackrel{v_2 \in V_2}{\operatorname{Tr}\mid_{\Gamma_2}v_2=0}} \mathcal J(u_1^{(n+1)} + v_2) &\\
u^{(n+1)}:=u_1^{(n+1)} + u_2^{(n+1)}\\
\tilde{u}_1^{(n+1)}:=\chi_1\cdot u^{(n+1)} \\
\tilde{u}_2^{(n+1)}:=\chi_2\cdot u^{(n+1)}.
\end{array}
\right.
\end{equation}
Note that we are minimizing over functions $v_i\in V_i$ for $i=1,2$ which vanish on the interior boundaries, i.e., $\operatorname{Tr}\mid_{\Gamma_i}v_i=0$. Moreover $u^{(n)}$ is the sum of the local minimizers $u_1^{(n)}$ and $u_2^{(n)}$, which are not uniquely determined on the overlapping part. Therefore we introduced a suitable correction by $\chi_1$ and $\chi_2$ in order to force the subminimizing sequences $(u_1^{(n)})_{n\in\mathbb{N}}$ and $(u_2^{(n)})_{n\in\mathbb{N}}$ to keep uniformly bounded. This issue will be explained in detail below, see Lemma \ref{bdtu12}. From the definition of $\chi_i$, $i=1,2$, it is clear that
$$
u_1^{(n+1)}+u_2^{(n+1)}=u^{(n+1)}=(\chi_1+\chi_2)u^{(n+1)}=\tilde{u}_1^{(n+1)} + \tilde{u}_2^{(n+1)}. 
$$
Note that in general $u_1^{(n)}\not=\tilde{u}_1^{(n)}$ and $u_2^{(n)}\not=\tilde{u}_2^{(n)}$. In \eqref{schw_sp} we use "$\approx$" (the approximation symbol) because in practice we never perform the exact minimization.
In the following section we discuss how to realize the approximation to the individual subspace minimizations.


\section{Local Minimization by Lagrange Multipliers}\label{localproblem}
Let us consider, for example, the subspace minimization on $\Omega_1$
\begin{equation}\label{funcsub1}
\operatorname{argmin}_{\stackrel{v_1 \in V_1}{\operatorname{Tr}\mid_{\Gamma_1}v_1=0}}  \mathcal  J(v_1 +u_2)=\operatorname{argmin}_{\stackrel{v_1 \in V_1}{\operatorname{Tr}\mid_{\Gamma_1}v_1=0}} \left\|Tv_1-(g-Tu_2)\right\|_2^2 + 2\alpha \left|\nabla(v_1+u_2 \mid_{\Omega_1})\right|(\Omega).
\end{equation}
First of all, observe that $\left\{u\in\mathcal{H}: \operatorname{Tr}\mid_{\Gamma_1}u=\operatorname{Tr}\mid_{\Gamma_1}u_2,\; \mathcal J(u)\leq C\right\}\subset \left\{\mathcal J\leq C\right\}$, hence the former set is also bounded by assumption (C) and the minimization problem \eqref{funcsub1} has solutions. 

It is useful to us to introduce an auxiliary functional $\mathcal J^s_1$ of $\mathcal{J}$, called the {\it surrogate functional} of $\mathcal J$ (cf. \cite{FS}): Assume $a,u_1 \in V_1$, $u_2 \in V_2$, and define
\begin{equation}
\label{surrfunc}
\mathcal J^s_1(u_1+ u_2, a) := \mathcal J(u_1+ u_2)+ \| u_1 -a\|_2^2 - \| T(u_1 -a)\|_2^2.
\end{equation}
A straightforward computation shows that
$$
\mathcal J^s_1(u_1+ u_2, a) = \| u_1 - (a + (T^*( g - T u_2 - T a))\mid_{\Omega_1})\|_2^2 + 2 \alpha \left|\nabla(u_1 + u_2)\right|(\Omega) + \Phi(a,g,u_2),
$$
where $\Phi$ is a function of $a,g,u_2$ only. Note that now the variable $u_1$ is not anymore effected by the action of $T$. {\nnew Consequently,} we want to realize an approximate solution to \eqref{funcsub1} by using the following algorithm: For $u_1^{(0)}=\tilde{u}_1^{(0)} \in V_1$,
\begin{equation}
\label{m2}
u_1^{(\ell+1)} = \operatorname{argmin}_{\stackrel{u_1 \in V_1}{\operatorname{Tr}\mid_{\Gamma_1}u_1=0}}  \mathcal  J^s_1(u_1 +u_2, u_1^{(\ell)}), \quad \ell \geq 0.
\end{equation}
Additionally in \eqref{m2} we can restrict the total variation on $\Omega_1$ only, since we have 
\begin{equation}\label{tvdecomp}
\begin{array}{lll}
\left|\nabla(u_1+u_2)\right|(\Omega) & = & \left|\nabla(u_1+u_2)\mid_{\Omega_1}\right|(\Omega_1) + c_2(u_2|_{(\Omega_2 \setminus \Omega_1) \cup \Gamma_1}).
\end{array}
\end{equation}
where we used \eqref{splitdom} and {\nnew the assumption} that $u_1$ vanishes on the interior boundary $\Gamma_1$. Hence \eqref{m2} is equivalent to
$$
\operatorname{argmin}_{\stackrel{u_1 \in V_1}{\operatorname{Tr}\mid_{\Gamma_1}u_1=0}}  \mathcal  J^s_1(u_1 +u_2, u_1^{(\ell)}) = \operatorname{argmin}_{\stackrel{u_1 \in V_1}{\operatorname{Tr}\mid_{\Gamma_1}u_1=0}} \left\|u_1-z_1\right\|_2^2 + 2\alpha \left| \nabla(u_1+u_2)\mid_{\Omega_1}\right|(\Omega_1),
$$
where $z_1= u_1^{(\ell)}+ (T^*( g - T u_2 - T u_1^{(\ell)}))\mid_{\Omega_1}$. Similarly the same arguments work for the second subproblem.

Before proving the convergence of this algorithm, we need to clarify first how to practically compute $u_1^{(\ell+1)}$ for $\tilde u_1^{(\ell)}$ given. To this end we need to introduce further notions and to recall some useful results.

\subsection{Generalized Lagrange multipliers for nonsmooth objective functions}
Let us begin this subsection with the notion of a {\it subdifferential} in finite dimensions.
\begin{definition}
\label{subdiff}
For a finite locally convex space $V$ and for a convex function $F:V \to \mathbb R \cup \{-\infty,+\infty\}$, we define the \emph{subdifferential} of $F$ at $x \in V$, as the set valued function
$$
\partial F(x) :=\begin{cases}
 \emptyset & \textrm{if } F(x) = \infty\\ 
 \{x^* \in V : \langle x^*,y-x\rangle + F(x) \leq F(y) \quad \forall y \in V\} & \textrm{otherwise}.
\end{cases}
$$
It is obvious from this definition that $0 \in \partial F(x)$ if and only if $x$ is a minimizer
of $F$. 
Since we deal with several spaces, namely, $\mathcal{H},V_i$, it will turn out to be useful to {\nnew sometimes distinguish} in which space the subdifferential is defined by imposing a subscript $\partial_V F$ for the subdifferential considered on the space $V$.
\end{definition}

We consider the following problem
\begin{equation}
\label{pbbp}
\operatorname{argmin}_{x \in V} \{ F(x): Gx=b\},
\end{equation}
where $ G:V \to V$ is a linear operator on $V$. We have the following useful result.

\begin{theorem}\cite[Theorem 2.1.4, p. 305]{HiLe}\label{huth}
Let $N=\left\{G^*\lambda:\;\lambda\in V\right\}=\operatorname{Range}(G^*)$. Then, $x_0\in\left\{x\in V:\; G(x)=b\right\}$ solves the constrained minimization problem \eqref{pbbp} if and only if
$$
0\in \partial F(x_0) + N.
$$
\end{theorem}


\subsection{Oblique thresholding}

We want to exploit Theorem \ref{huth} in order to produce an algorithmic solution to each iteration step (\ref{m2}), which practically stems from the solution of a problem of this type
$$
\operatorname{argmin}_{\stackrel{u_1 \in V_1}{\operatorname{Tr}\mid_{\Gamma_1}u_1=0}}  \| u_1 -z_1\|_{2}^2 + 2 \alpha \left|\nabla(u_1+u_2 \mid_{\Omega_1})\right|(\Omega_1).
$$
It is well-known how to solve this problem if {\nnew $u_2\equiv 0$ in $\bar{\Omega}_1$} and the trace condition is not imposed. {\nnew For the general case we propose the following solution strategy.} In what follows all the involved quantities are restricted to $\Omega_1$, e.g., $u_1=u_1\mid_{\Omega_1}, u_2=u_2\mid_{\Omega_1}$.

\begin{theorem}[Oblique thresholding]
\label{main1}
For $u_2 \in V_2$ and for $z_1 \in V_1$ the following statements are equivalent:
\begin{itemize}
\item[(i)] $u_1^* = \operatorname{argmin}_{\stackrel{u_1 \in V_1}{\operatorname{Tr}\mid_{\Gamma_1}u_1=0}}  \| u_1 -z_1\|_{2}^2 + 2 \alpha \left|\nabla(u_1+u_2)\right|(\Omega_1)$;
\item[(ii)] there exists $\eta \in \operatorname{Range} (\operatorname{Tr}\mid_{\Gamma_1})^* = \left\{\eta\in V_1 \textrm{ with } \operatorname{supp}(\eta)=\Gamma_1\right\}$ such that
$0 \in u^*_1 -(z_1- \eta) + \alpha \partial_{V_1} \left|\nabla(\cdot + u_2)\right|(\Omega_1)(u_1^*)$;
\item[(iii)] there exists $\eta \in V_1$ with $\operatorname{supp}(\eta)=\Gamma_1$ such that $u_1^* = (I- P_{\alpha K})(z_1 + u_2 -\eta)-u_2 \in V_1$ and $\operatorname{Tr}\mid_{\Gamma_1}u_1^*=0$;
\item[(iv)] there exists $\eta \in V_1$  with $\operatorname{supp}(\eta)=\Gamma_1$ such that  $\operatorname{Tr}\mid_{\Gamma_1} \eta = \operatorname{Tr}\mid_{\Gamma_1} z_1 + \operatorname{Tr}\mid_{\Gamma_1} P_{\alpha K} (\eta - (z_1 + u_2))$ or equivalently $\eta = (\operatorname{Tr}\mid_{\Gamma_1})^* \operatorname{Tr}\mid_{\Gamma_1} \left(z_1 + P_{\alpha K} (\eta - (z_1 + u_2))\right)$.
\end{itemize}
\end{theorem}

We call the solution operation provided by this theorem an {\it oblique thresholding}, in analogy to the terminology in \cite{dadede04}, because it performs a thresholding of the derivatives, i.e., it sets to zero most of the derivatives of $u=u_1 + u_2 \approx z_1$ on $\Omega_1$, provided $u_2$ which is a fixed vector in $V_2$.

\begin{proof}
Let us show the equivalence between (i) and (ii). The problem in (i) can be reformulated as
\begin{equation}
u_1^* = \operatorname{argmin}_{u_1 \in V_1} \{ F(u_1):=\| u_1 -z_1\|_{2}^2 + 2 \alpha \left|\nabla(u_1+u_2)\right|(\Omega_1), \operatorname{Tr}\mid_{\Gamma_1}u_1=0 \}.
\end{equation}

Recall that $\operatorname{Tr}\mid_{\Gamma_1}: V_1 \rightarrow \mathbb{R}^{\Gamma_1}$ is a surjective map with closed range. This means that $(\operatorname{Tr}\mid_{\Gamma_1})^*$ is injective and that $\operatorname{Range} (\operatorname{Tr}\mid_{\Gamma_1})^* = \left\{\eta\in V_1 \textrm{ with } \operatorname{supp}(\eta)=\Gamma_1\right\}$ is closed. Using Theorem \ref{huth} the optimality of $u_1^*$ is equivalent to the existence of $\eta \in \operatorname{Range}(\operatorname{Tr}\mid_{\Gamma_1})^*$ such that 
\begin{equation}
0\in \partial_{V_1}F(u_1^*) + 2 \eta.
\end{equation}
Due to the continuity of $\|u_1-z_1\|_{2}^2$ in $V_1$, we have, by \cite[Proposition 5.6]{ET}, that 
\begin{equation}
\partial_{V_1} F(u_1^*) = 2 (u_1^*-z_1) + 2 \alpha \partial_{V_1}\left|\nabla(\cdot+u_2)\right|(\Omega_1)(u_1^*).
\end{equation}
Thus, the optimality of $u_1^*$ is equivalent to
\begin{equation}
0 \in u_1^* - z_1 + \eta + \alpha \partial_{V_1} \left|\nabla(\cdot+u_2)\right|(\Omega_1)(u_1^*).
\end{equation}
This concludes the equivalence of (i) and (ii).
Let us show now that (iii) is equivalent to (ii). 
The condition in (iii) can be rewritten as
$$
\xi^* = (I- P_{\alpha K})(z_1+u_2 - \eta), \hspace{1cm} \xi^*=u_1^*+u_2. 
$$
Since $\left|\nabla(\cdot)\right| \geq 0 $ is 1-homogeneous and lower-semicontinuous, by \cite[Example 4.2.2]{FS}, the latter is equivalent to
$$
0 \in  \xi^* - (z_1+u_2 - \eta) + \alpha \partial_{V_1} \left|\nabla(\cdot)\right|(\Omega_1)(\xi^*),
$$
and equivalent to (ii). Note that in particular we have $\partial_{V_1} \left|\nabla(\cdot)\right|(\Omega_1)(\xi^*)=\partial_{V_1} \left|\nabla(\cdot+u_2)\right|(\Omega_1)(u_1^*)$, which is easily shown by a direct computation from the definition of subdifferential. We prove now the equivalence between (iii) and (iv). We have
\begin{eqnarray*}
u_1^* &=& (I- P_{\alpha K})(z_1 + u_2 -\eta)-u_2 \in V_1, \quad \eta\in V_1 \text{ with } \supp(\eta) = \Gamma_1,\operatorname{Tr}\mid_{\Gamma_1}u_1^*=0 \\
&=& z_1 - \eta - P_{\alpha K}(z_1 + u_2 -\eta).
\end{eqnarray*}
By applying $\operatorname{Tr}\mid_{\Gamma_1}$ to both sides of the latter equality we get
$$
0 = \operatorname{Tr}\mid_{\Gamma_1} z_1 - \operatorname{Tr}\mid_{\Gamma_1} \eta - \operatorname{Tr}\mid_{\Gamma_1} P_{\alpha K}(z_1 + u_2 -\eta).
$$
By observing that {\nnew $-\operatorname{Tr}\mid_{\Gamma_1} P_{\alpha K}(\xi)=\operatorname{Tr}\mid_{\Gamma_1} P_{\alpha K}(-\xi)$,} we obtain the fixed point equation
\begin{equation}
\label{fixpt}
\operatorname{Tr}\mid_{\Gamma_1} \eta = \operatorname{Tr}\mid_{\Gamma_1} z_1 + \operatorname{Tr}\mid_{\Gamma_1} P_{\alpha K} (\eta - (z_1 + u_2)).
\end{equation}
Conversely, since all the considered quantities in 
$$
(I- P_{\alpha K})(z_1 + u_2 -\eta) - u_2
$$
are in $V_1$, the whole expression is an element in $V_1$ and hence $u_1^*$ as defined in (iii) is an element in $V_1$ and $\operatorname{Tr}\mid_{\Gamma_1}u_1^*=0$. {\nnew This shows the equivalence between (iii) and (iv) and therewith finishes the proof.}
\end{proof}

We wonder now whether any of {\nnew the conditions in Theorem \ref{main1}} is indeed practically satisfied. {\nnew In particular, we} want to show that $\eta \in V_1$ as in (iii) or (iv) of the previous theorem is provided as the limit of the following iterative algorithm:
\begin{equation}
\label{fixptit1}
\eta^{(0)} \in V_1,\supp \eta^{(0)}=\Gamma_1 \quad \eta^{(m+1)} =(\operatorname{Tr}\mid_{\Gamma_1})^*  \operatorname{Tr}\mid_{\Gamma_1}\left(z_1 + P_{\alpha K}(\eta^{(m)}- (z_1 + u_2))\right), \quad m \geq 0.
\end{equation}

\begin{proposition}\label{fixptitpr}
The following statements are equivalent:
\begin{itemize}
\item[(i)] there exists $\eta \in V_1$ such that
$\eta = (\operatorname{Tr}\mid_{\Gamma_1})^*\operatorname{Tr}\mid_{\Gamma_1}\left(z_1 + P_{\alpha K} (\eta - (z_1 + u_2))\right)$ (which is in turn the condition (iv) of Theorem \ref{main1})
\item[(ii)] the iteration (\ref{fixptit1}) converges to any $\eta \in V_1$ that satisfies (\ref{fixpt}).
\end{itemize}
\end{proposition}

For the proof of this Proposition we need to recall some well-known notions and results.
\begin{definition}
A nonexpansive map $\mathcal T:\mathcal H \to \mathcal H$ is strongly nonexpansive if for $(u_n-v_n)_n$ bounded and $\| \mathcal T(u_n) - \mathcal T(v_n)\|_{2} - \| u_n-v_n\|_{2} \to 0$ we have
$$
u_n-v_n - (\mathcal T(u_n) - \mathcal T(v_n)) \to 0, \quad n \to \infty.
$$
\end{definition}
\begin{proposition}[Corollaries 1.3, 1.4, and 1.5 \cite{br}]
\label{bbl}
Let $\mathcal T:\mathcal H \to \mathcal H$ be a strongly nonexpansive map. Then $\operatorname{fix} \mathcal T =\{ u \in \mathcal H: \mathcal T(u)=u\} \neq \emptyset$ if and only if $(\mathcal T^n u)_n$ converges to a fixed point $u_0 \in \operatorname{fix} \mathcal T$ for any choice of $u \in \mathcal H$.
\end{proposition}

\begin{proof}(Proposition \ref{fixptitpr})
Projections onto convex sets are strongly nonexpansive \cite[Corollary 4.2.3]{bbl}. Moreover, the composition of strongly nonexpansive maps is strongly nonexpansive \cite[Lemma 2.1]{br}. By an application of Proposition \ref{bbl} we immediately have the result, since any map of the type $\mathcal T(\xi) = Q(\xi) + \xi_0$ is strongly nonexpansive whenever $Q$ is (this is a simple observation from the definition of strongly nonexpansive maps).
Indeed, we are looking for fixed points of $\eta = (\operatorname{Tr}\mid_{\Gamma_1})^* \operatorname{Tr}\mid_{\Gamma_1}(z_1 + P_{\alpha K} (\eta - (z_1 + u_2)))$ or, equivalently, of $\xi = \underbrace{(\operatorname{Tr}\mid_{\Gamma_1} )^*\operatorname{Tr}\mid_{\Gamma_1} P_{\alpha K}}_{:=Q} (\xi) - \underbrace{((\operatorname{Tr}\mid_{\Gamma_1} )^*\operatorname{Tr}\mid_{\Gamma_1} u_2)}_{:=\xi_0}$, {\nnew where $\xi = (\operatorname{Tr}\mid_{\Gamma_1})^* \operatorname{Tr}\mid_{\Gamma_1} (\eta - (z_1 + u_2))$}.
\end{proof}

\subsection{Convergence of the subspace minimization}
From the results of the previous section it follows that the iteration \eqref{m2} can be explicitly computed by
\begin{equation}\label{m3}
u_1^{(\ell+1)} = S_{\alpha}(u_1^{(\ell)} + T^* (g - T u_2 - T   u_1^{(\ell)}) + u_2 - \eta^{(\ell)}) - u_2,
\end{equation} 
where $S_{\alpha}:=I-P_{\alpha K}$ and $\eta^{(\ell)}\in V_1$ is any solution of the fixed point equation
$$
\eta = (\operatorname{Tr}\mid_{\Gamma_1})^*\operatorname{Tr}\mid_{\Gamma_1} \left((u_1^{(\ell)} + T^* (g - T u_2 - T   u_1^{(\ell)})) - P_{\alpha K}(u_1^{(\ell)} + T^* (g - T u_2 - T   u_1^{(\ell)}+ u_2 -\eta))\right).
$$
The computation of $\eta^{(\ell)}$ can be implemented by the algorithm \eqref{fixptit1}.

\begin{proposition}\label{main2}
Assume $u_2 \in V_2$ and $\| T\| <1$. Then the iteration \eqref{m3} converges to a solution $u^*_1 \in V_1$ of \eqref{funcsub1} for any initial choice of $u_1^{(0)} \in V_1$.
\end{proposition}

The proof of this proposition is standard, see \cite{CW,dadede04,FS}.

Let us conclude this section mentioning that all the results presented here hold symmetrically for the minimization on $V_2$, and that the notations should be just adjusted accordingly.


\section{Convergence of the Sequential Alternating Subspace Minimization}\label{Convsec}
In this section we want to prove the convergence of the algorithm (\ref{schw_sp}) to minimizers of $\mathcal{J}$. In order to do that, we need a characterization of solutions of the minimization problem \eqref{functotalfin} as the one provided in \cite[Proposition 4.1]{Ve01} for the continuous setting. We specify the arguments in \cite[Proposition 4.1]{Ve01} for our discrete setting and we highlight the significant differences with respect to the continuous one.


\subsection{Characterization of Solutions}

We make the following assumptions:
\begin{itemize}
\item[$(A_\varphi)$] $\varphi: \R \to \R$ is a convex function, nondecreasing in $\R^+$ such that
	\begin{itemize}
	\item[(i)] $\varphi(0)=0$.
	\item[(ii)] There exist $c>0$ and $b\geq 0$ such that $cz-b\leq \varphi(z)\leq cz+b,$ for all $z\in \R^+$.
	\end{itemize}
\end{itemize}
The particular example we have in mind is simply $\varphi(s)=s$, but we keep a more general notation for uniformity with respect to the continuous version in \cite[Proposition 4.1]{Ve01}.
In this section we are concerned with the following more general minimization problem
\begin{equation}\label{problem}
\operatorname{argmin}_{u\in\mathcal{H}}\{\mathcal{J}_{\varphi}(u):=\|Tu-g\|_2^2 + 2\alpha \varphi(|\nabla u|)(\Omega)\}
\end{equation}
where $g\in\mathcal{H}$ is a datum, $\alpha > 0$ is a fixed constant (in particular for  $\varphi(s)=s$).

To characterize the solution of the minimization problem \eqref{problem} we use duality results from \cite{ET}. Therefore we recall the definition of the {\it conjugate (or Legendre transform)} of a function (for example see \cite[Def. 4.1, pag. 17]{ET}):
\begin{definition}\label{Def.conjugate}
Let $V$ and $V^*$ be two vector spaces placed in the duality by a bilinear pairing denoted by $\langle\cdot,\cdot\rangle$ and $\phi:V\to\R$ be a convex function. The \emph{conjugate function (or Legendre transform)} $\phi^*:V^*\to\R$ is defined by
$$
\phi^*(u^*) = \sup_{u\in V}\{\langle u,u^* \rangle - \phi(u)\}.
$$
\end{definition}

\begin{proposition}\label{finitVese}
Let {\nnew $\zeta,u \in \mathcal{H}$}. If the assumption $(A_\varphi)$ is fulfilled, then $\zeta \in \partial \mathcal{J}_{\varphi}(u)$ if and only if there exists $M=(M_0,\bar{M})\in\mathcal{H}\times\mathcal{H}^d$, $\frac{|\bar M(x)|}{2\alpha}\leq c_1\in[0,+\infty)$ for all $x\in\Omega$ such that
\begin{eqnarray}
\langle\bar M(x),(\nabla u)(x) \rangle_{\R^d} + 2\alpha \varphi(|(\nabla u)(x)|) + 2 \alpha \varphi_1^*\left(\frac{|\bar M(x)|}{2 \alpha}\right) &=& 0 \quad \text{for all } x\in\Omega \label{cond1}\\
T^*M_0 - \operatorname{div} \bar M + {\nnew \zeta} &=& 0\label{cond2}\\
-M_0 = 2(Tu-g)\label{cond3},
\end{eqnarray}
where $\varphi_1^*$ is the conjugate function of $\varphi_1$ defined by $\varphi_1(s)=\varphi(|s|)$, for $s\in\R$.

If additionally $\varphi$ is differentiable and $|(\nabla u)(x)|\not=0$ for $x \in\Omega$, then we can compute $\bar M$ as
\begin{equation}\label{Mbar}
\bar{M}(x) = -2\alpha \frac{\varphi'(|(\nabla u)(x)|)}{|(\nabla u)(x)|}(\nabla u)(x).
\end{equation}
\end{proposition}

The proof of this proposition specifies the one of \cite[Proposition 4.1]{Ve01} to our discrete setting, it is technical, and it is deferred to the Appendix. 

\begin{remark}\label{remarkk}
\begin{enumerate}
\item[(i)] For $\varphi(s)=s$ the function $\varphi_1$ from Proposition \ref{finitVese} turns out to be $\varphi_1(s)=|s|$. Its conjugate function $\varphi_1^*$ is then given by 
$$
\varphi_1^*(s^*)=\sup_{s\in\R}\{\langle s^*,s \rangle - |s|\}=
\begin{cases}
0 \quad &\text{for }|s^*|\leq 1\\
\infty	\quad &\text{else}
\end{cases}.
$$
Hence condition \eqref{cond1} specifies as follows
$$
\langle\bar{M}(x),(\nabla u)(x) \rangle_{\R^d} + 2\alpha |(\nabla u)(x)| = 0
$$
and, directly from the proof of Proposition \ref{finitVese} in the Appendix, $|\bar{M}(x)|\leq 2\alpha$ for all $x\in\Omega$.
\item[(ii)] We want to highlight a few important differences with respect to the continuous case. Due to our definition of the gradient and its relationship with the divergence operator  $-\operatorname{div}=\nabla^*$ no boundary conditions are needed. Therefore condition (10) of \cite[Proposition 4.1]{Ve01} has no discrete correspondent in our setting. The continuous total variation of a function can be decomposed into an absolute continuous part with respect to the Lebesgue measure and a singular part, whereas no singular part appears in the discrete setting. Therefore condition (6) and (7) of \cite[Proposition 4.1]{Ve01} have not a discrete correspondent either.
\item[(iii)] An interesting consequence of Proposition \ref{finitVese} {\nnew is that} the map $S_\alpha =(I- P_{\alpha K})$ is bounded, i.e., $\| S_\alpha(z^k)\|_2 \to \infty$ {\nnew if and only if} $\| z^k \|_2 \to \infty$, for $k \to \infty$. In fact, since
\begin{eqnarray*}
S_\alpha (z) &=& \arg \min_{u \in \mathcal H} \| u - z \|_2^2 + 2 \alpha | \nabla u|(\Omega),
\end{eqnarray*}
from \eqref{cond2} and \eqref{cond3}, we immediately obtain
$$
S_\alpha (z) = z - \frac{1}{2} \operatorname{div} \bar M,
$$
and $\bar M$ is uniformly bounded.
\end{enumerate}
\end{remark}

\subsection{Convergence properties}
We return to the {\nnew sequential} algorithm (\ref{schw_sp}). Let us explicitly express the algorithm as follows:
Pick an initial $V_1 + V_2 \ni \tilde u_1^{(0)}+ \tilde u_2^{(0)} : = u^{(0)} \in \mathcal{H}$, for example, $\tilde u_i^{(0)}=0, i=1,2$, and iterate
\begin{equation}
\label{schw_sp:it2}
\left \{ 
\begin{array}{ll}
\left \{ 
\begin{array}{ll}
u_1^{(n+1,0)} = \tilde u_1^{(n)}&\\
u_1^{(n+1,\ell+1)} =  \operatorname{argmin}_{\stackrel{u_1 \in V_1}{\operatorname{Tr}\mid_{\Gamma_1}u_1=0}} \mathcal J_1^s(u_1+ \tilde u_2^{(n)}, u_1^{(n+1,\ell)}) & \ell=0,\dots, L-1\\
\end{array}\right. &\\
\left \{ 
\begin{array}{ll}
u_2^{(n+1,0)} = \tilde{u}_2^{(n)}&\\
u_2^{(n+1,m+1)} = \operatorname{argmin}_{\stackrel{u_2 \in V_2}{\operatorname{Tr}\mid_{\Gamma_2}u_2=0}} \mathcal J_2^s(u_1^{(n+1,L)}+ u_2, u_2^{(n+1,m)}) &m=0,\dots, M -1\\
\end{array}\right. &\\
u^{(n+1)}:=u_{1}^{(n+1,L)} + u_{2}^{(n+1,M)}\\
\tilde{u}_1^{(n+1)}:=\chi_1\cdot u^{(n+1)} \\
\tilde{u}_2^{(n+1)}:=\chi_2\cdot u^{(n+1)}.
\end{array}
\right.
\end{equation} 
Note that we do prescribe a finite number $L$ and $M$ of inner iterations for each subspace respectively and that $u^{(n+1)}=\tilde u_1^{(n+1)} + \tilde u_2^{(n+1)}$, with $u_i^{(n+1)}\not=\tilde u_i^{(n+1)}$, $i=1,2,$ in general. In this section we want to prove its convergence for any choice of $L$ and $M$.\\

Observe that, for $a \in V_i$ and $\|T\|<1$, 
\begin{equation}
\| u_i- a\|_{2}^2 - \| Tu_i- Ta\|^2_{2} \geq C \| u_i- a\|_{2}^2,
\end{equation}
for $C=(1-\|T\|^2)>0$. Hence
\begin{equation}
\label{decr}
\mathcal {J}(u) = \mathcal {J}_i^s(u,u_i) \leq \mathcal {J}_i^s(u,a), 
\end{equation}
and
\begin{equation}
\label{decr2}
\mathcal {J}_i^s(u,a) - \mathcal {J}_i^s(u,u_i) \geq C \|u_i-a\|^2_2.
\end{equation}

\begin{proposition}[Convergence properties]
\label{weak-conv}
Let us assume that $\|T\|<1$. The algorithm in (\ref{schw_sp:it2}) produces a sequence $(u^{(n)})_{n\in \mathbb{N}}$ in  $\mathcal{H}$ with the following properties:
\begin{itemize}
\item[(i)] $\mathcal{J}(u^{(n)}) > \mathcal{J}(u^{(n+1)})$ for all $n \in \mathbb{N}$ (unless $u^{(n)}= u^{(n+1)}$);
\item[(ii)] $\lim_{n \to \infty} \| u^{(n+1)} -  u^{(n)}\|_{2} = 0$;
\item[(iii)] the sequence  $(u^{(n)})_{n\in \mathbb{N}}$ has subsequences which converge in $\mathcal{H}$.
\end{itemize}
\end{proposition}

\begin{proof}
Let us first observe that 
\begin{eqnarray*}
 \mathcal {J}(u^{(n)})= \mathcal {J}_1^s(\tilde u^{(n)}_1+ \tilde u^{(n)}_2, \tilde u^{(n)}_1)&=& \mathcal {J}^s_1(\tilde u_1^{(n)}  + \tilde u_2^{(n)}, u_1^{(n+1,0)}).
\end{eqnarray*}
By definition of $u_1^{(n+1,1)}$ and the minimal properties of $u_1^{(n+1,1)}$  in (\ref{schw_sp:it2}) we have
$$
\mathcal {J}_1^s(\tilde u_1^{(n)}  + \tilde u_2^{(n)}, u_1^{(n+1,0)}) \geq \mathcal {J}^s_1(u_1^{(n+1,1)}   + \tilde u_2^{(n)}, u_1^{(n+1,0)}).
$$
From (\ref{decr}) we have
$$
\mathcal {J}^s_1(u_1^{(n+1,1)} + \tilde u_2^{(n)}, u_1^{(n+1,0)}) \geq\mathcal {J}^s_1(u_1^{(n+1,1)} + \tilde u_2^{(n)}, u_1^{(n+1,1)}) = \mathcal{J}(u_1^{(n+1,1)} + \tilde u_2^{(n)}).
$$
Putting in line these inequalities we obtain 
$$
\mathcal {J}(u^{(n)}) \geq \mathcal {J}(u_1^{(n+1,1)} + \tilde u_2^{(n)}).
$$
In particular, from (\ref{decr2}) we have
$$
\mathcal {J}(u^{(n)}) - \mathcal {J}(u_1^{(n+1,1)} + \tilde u_2^{(n)}) \geq C \| u_1^{(n+1,1)} - u_1^{(n+1,0)}\|_{2}^2.
$$
After $L$ steps we conclude the estimate
\begin{eqnarray*}
\mathcal {J}(u^{(n)}) &\geq& \mathcal {J} (u_1^{(n+1,L)} + \tilde u_2^{(n)}),
\end{eqnarray*}
and
$$
\mathcal {J}(u^{(n)}) -  \mathcal {J}(u_1^{(n+1,L)} + \tilde u_2^{(n)})\geq C \sum_{\ell=0}^{L-1}\| u_1^{(n+1,\ell+1)} - u_1^{(n+1,\ell)}\|_{2}^2.
$$
By definition of $u_2^{(n+1,1)}$ and its minimal properties we have
\begin{eqnarray*}
&& \mathcal {J}(u_1^{(n+1,L)} + \tilde u_2^{(n)})\geq  \mathcal {J}^s_2(u_1^{(n+1,L)} + u_2^{(n+1,1)}, u_2^{(n+1,0)} ).
\end{eqnarray*}
By similar arguments as above we finally find the decreasing estimate
\begin{equation}
\label{decr3}
\mathcal {J}(u^{(n)}) \geq  \mathcal {J}(u_1^{(n+1,L)} + u_2^{(n+1,M)}) = \mathcal {J}(u^{(n+1)})=\mathcal{J}(\tilde u_1^{(n+1)}+\tilde u_2^{(n+1)}),
\end{equation}
and
$$
\mathcal {J}(u^{(n)}) - \mathcal {J}(u^{(n+1)}) 
$$
\begin{equation}
\label{coercive}\geq C \left ( \sum_{\ell=0}^{L-1}\| u_1^{(n+1,\ell+1)} - u_1^{(n+1,\ell)}\|_{2}^2 + \sum_{m=0}^{M-1}\| u_2^{(n+1,m+1)} - u_2^{(n+1,m)}\|_{2}^2 \right ),
\end{equation}
which verifies $(i)$.\\
From (\ref{decr3}) we have $\mathcal {J}(u^{(0)}) \geq \mathcal {J}(u^{(n)})$. By the coerciveness condition (C) $(u^{(n)})_{n \in \mathbb{N}}$ is uniformly bounded in $\mathcal{H}$, hence there exists a convergent subsequence $(u^{(n_k)})_{k \in \mathbb{N}}$ and hence $(iii)$ holds. 
Let us denote $u^{(\infty)}$ the limit of the subsequence. For simplicity, we rename such a subsequence by $(u^{(n)})_{n \in \mathbb{N}}$. Moreover, since the sequence $(\mathcal {J}(u^{(n)}))_{n \in \mathbb{N}}$ is monotonically decreasing and bounded from below by 0, it is also convergent. From (\ref{coercive}) and the latter convergence we deduce 

\begin{equation}
\label{asymp_reg}
\left ( \sum_{\ell=0}^{L-1}\| u_1^{(n+1,\ell+1)} - u_1^{(n+1,\ell)}\|_{2}^2 + \sum_{m=0}^{M-1}\| u_2^{(n+1,m+1)} - u_2^{(n+1,m)}\|_{2}^2 \right ) \rightarrow 0, \quad n \to \infty.
\end{equation}
In particular, by the standard inequality $(a^2+b^2) \geq \frac{1}{2} (a+b)^2$ for $a,b>0$ and the triangle inequality, we have also
\begin{equation}
\label{asymp_reg2}
\| u^{(n)} -  u^{(n+1)} \|_{2} \rightarrow 0, \quad n\to \infty.
\end{equation}
This gives $(ii)$ and completes the proof.
\end{proof}
The use of the partition of unity $\{\chi_1, \chi_2\}$ allows not only to guarantee the boundedness of $(u^{(n)})_{n \in \mathbb N}$, but also of the sequences $(\tilde u_1^{(n)})_{n\in \mathbb{N}}$ and $(\tilde u_2^{(n)})_{n\in \mathbb{N}}$.
\begin{lemma}\label{bdtu12}
The sequences $(\tilde u_1^{(n)})_{n\in \mathbb{N}}$ and $(\tilde u_2^{(n)})_{n\in \mathbb{N}}$ produced by the algorithm \eqref{schw_sp:it2} are bounded, i.e., there exists a constant $\tilde C >0$ such that $\|\tilde u_i^{(n)}\|_2\leq \tilde C$ for $i=1,2$.
\end{lemma}
\begin{proof}
From the boundedness of $(u^{(n)})_{n \in \mathbb N}$ we have 
$$
\|\tilde u_i^{(n)}\|_2 = \|\chi_i u^{(n)}\|_2 \leq \kappa \|u^{(n)}\|_2 \leq \tilde C \hspace{0.5cm} \text{for } i=1,2.
$$
\end{proof}
From Remark \ref{remarkk} (iii) we can also show the following auxiliary lemma.
\begin{lemma}\label{eta_bd}
The sequences {\nnew $(\eta_1^{(n,L)})_{n}$ and $^(\eta_2^{(n,M)})_{n}$} are bounded.
\end{lemma}
\begin{proof}
From previous considerations we know that
\begin{align*}
u_1^{(n,L)} &= S_\alpha(z_1^{(n,L-1)} + \tilde u_2^{(n-1)} - \eta_1^{(n,L)}) - \tilde u_2^{(n-1)}\\
u_2^{(n,M)} &= S_\alpha(z_2^{(n,M-1)} + u_1^{(n,L)} - \eta_2^{(n,M)}) - u_1^{(n,L)}.
\end{align*} 
Assume $(\eta_1^{(n,L)})_n$ were unbounded, then by Remark \ref{remarkk} (iii), also $S_\alpha(z_1^{(n,L-1)} + \tilde u_2^{(n-1)} - \eta_1^{(n,L)})$ would be unbounded. Since $(\tilde u_2^{(n)})_n$ and $(u_1^{(n,L)})_n$ are bounded by Lemma \ref{bdtu12} and formula \eqref{asymp_reg}, we have a contradiction. Thus $(\eta_1^{(n,L)})_n$ has to be bounded. With the same argument we can show that $(\eta_2^{(n,M)})_n$ is bounded.
\end{proof}

We can eventually show the convergence of the algorithm to minimizers of $\mathcal J$.

\begin{theorem}[Convergence to minimizers]\label{convergence:proof}
Assume $\|T\|<1$. Then accumulation points of the sequence $(u^{(n)})_{n\in \mathbb{N}}$ produced by algorithm (\ref{schw_sp:it2}) are minimizers of $\mathcal{J}$. If $\mathcal{J}$ has a unique minimizer then the sequence $(u^{(n)})_{n\in \mathbb{N}}$ converges to it.
\end{theorem}

\begin{proof}
Let us denote $u^{(\infty)}$ the limit of a subsequence. For simplicity, we rename such a subsequence by $(u^{(n)})_{n \in \mathbb{N}}$. From Lemma \ref{bdtu12} we know that $(\tilde{u}_1^{(n)})_{n\in\mathbb{N}}$, $(\tilde{u}_2^{(n)})_{n\in\mathbb{N}}$ and consequently $(u_1^{(n,L)})_{n\in\mathbb{N}}$,$(u_2^{(n,M)})_{n\in\mathbb{N}}$ are bounded. So the limit $u^{(\infty)}$ can be written as
\begin{equation}\label{limit}
u^{(\infty)} = u_1^{(\infty)} + u_2^{(\infty)} = \tilde{u}_1^{(\infty)} + \tilde{u}_2^{(\infty)}
\end{equation}
where $u_1^{(\infty)}$ is the limit of $(u_1^{(n,L)})_{n\in\mathbb{N}}$, $u_2^{(\infty)}$ is the limit of $(u_2^{(n,M)})_{n\in\mathbb{N}}$, and $\tilde{u}_i^{(\infty)}$ is the limit of $(\tilde{u}_i^{(n)})_{n\in\mathbb{N}}$ for $i=1,2$. Now we show that $\tilde{u}_2^{(\infty)} = u_2^{(\infty)}$.
{\nnew By using the triangle inequality, from \eqref{asymp_reg} it directly follows that}
\begin{equation}\label{A1}
\|u_2^{(n+1,M)}-\tilde{u}_2^{(n)}\|_2 \to 0, \quad n\to \infty.
\end{equation}
Moreover, since $\chi_2\in V_2$ is a fixed vector which is independent of $n$, we obtain from Proposition \ref{weak-conv} $(ii)$ that
$$
\| \chi_2(u^{(n)} -  u^{(n+1)}) \|_{2} \rightarrow 0, \quad n\to \infty,
$$
and hence
\begin{equation}\label{A2}
\| \tilde{u}_2^{(n)} -  \tilde{u}_2^{(n+1)} \|_{2} \rightarrow 0, \quad n\to \infty.
\end{equation}
Putting \eqref{A1} and \eqref{A2} together and noting that
$$
\|u_2^{(n+1,M)}-\tilde{u}_2^{(n)}\|_2 + \| \tilde{u}_2^{(n)} -  \tilde{u}_2^{(n+1)} \|_{2} \geq \| u_2^{(n+1,M)} -  \tilde{u}_2^{(n+1)} \|_{2}
$$
we have
\begin{equation}\label{samelimit}
\| u_2^{(n+1,M)} -  \tilde{u}_2^{(n+1)} \|_{2} \rightarrow 0, \quad n\to \infty,
\end{equation}
which means that the sequences $(u_2^{(n,M)})_{n\in\mathbb{N}}$ and $(\tilde{u}_2^{(n)})_{n\in\mathbb{N}}$ have the same limit, i.e., $\tilde{u}_2^{(\infty)} = u_2^{(\infty)}$, which we denote by $u_2^{(\infty)}$. Then from \eqref{samelimit} and \eqref{limit} it directly follows that $\tilde{u}_1^{(\infty)}=u_1^{(\infty)}$.

As in the proof of the oblique thresholding theorem we set
$$
F_1(u_1^{(n+1,L)}):=\|u_1^{(n+1,L)} - z_1^{(n+1,L)}\|_2^2 + 2 \alpha |\nabla(u_1^{(n+1,L)}+ \tilde u_2^{(n)}\Big{|}_{\Omega_1})|(\Omega_1)
$$
where
$$
{\nnew z_1^{(n+1,L)}}:= u_1^{(n+1,L-1)} + (T^* ( g - T  \tilde{u}_2^{(n)} - T u_1^{(n+1,L-1)}))\Big{|}_{\Omega_1}.
$$
The optimality condition for $u_1^{(n+1,L)}$ is
$$
0\in\partial_{V_1} F_1(u_1^{(n+1,L)}) + 2 \eta_1^{(n+1,L)}
$$
where
\begin{equation*}
{\nnew \eta_1^{(n+1,L)}} = (\operatorname{Tr}\mid_{\Gamma_1})^* \operatorname{Tr}\mid_{\Gamma_1}\left((z_1^{(n+1,L)}) + P_{\alpha K}(\eta_1^{(n+1,L)}-z_1^{(n+1,L)} - \tilde{u}_2^{(n)})\right).
\end{equation*}

In order to use the characterization of elements in the subdifferential of $|\nabla u|(\Omega)$, i.e., Proposition \ref{finitVese}, we have to rewrite the minimization problem for $F_1$. More precisely, we define
$$
\hat F_1(\xi_1^{(n+1,L)}):=\|\xi_1^{(n+1,L)} - \tilde u_2^{(n)}\Big{|}_{\Omega_1} - z_1^{(n+1,L)}\|_2^2 + 2 \alpha |\nabla(\xi_1^{(n+1,L)})|(\Omega_1)
$$
for $\xi_1^{(n+1,L)} \in V_1$ with $\operatorname{Tr}\mid_{\Gamma_1}\xi_1^{(n+1,L)} = \tilde u_2^{(n)}$. Then the optimality condition for $\xi_1^{(n+1,L)}$ is
\begin{equation}\label{O1}
0\in\partial \hat F_1(\xi_1^{(n+1,L)}) + 2 \eta_1^{(n+1,L)}
\end{equation}
Note that indeed $\xi_1^{(n+1,L)}$ is optimal if and only if $u_1^{(n+1,L)}=\xi_1^{(n+1,L)}-\tilde u_2^{(n)}\Big{|}_{\Omega_1}$ is optimal.

Analogously we define
$$
\hat F_2(\xi_2^{(n+1,M)}):=\|\xi_2^{(n+1,M)} - u_1^{(n+1,L)}\Big{|}_{\Omega_2} - z_2^{(n+1)}\|_2^2 + 2 \alpha |\nabla(\xi_2^{(n+1,M)})|(\Omega_2)
$$
for $\xi_2^{(n+1,M)} \in V_2$ with $\operatorname{Tr}\mid_{\Gamma_2}\xi_2^{(n+1,M)} = u_1^{(n+1,L)}$, and the optimality condition for $\xi_2^{(n+1,M)}$ is
\begin{equation}\label{O2}
0\in\partial \hat F_2(\xi_2^{(n+1,M)}) + 2 \eta_2^{(n+1,M)}
\end{equation}
where
\begin{equation*}
{\nnew \eta_2^{(n+1,M)}} = (\operatorname{Tr}\mid_{\Gamma_2})^*\operatorname{Tr}\mid_{\Gamma_2}\left((z_2^{(n+1,M)}) + P_{\alpha K}(\eta_2^{(n+1,M)}-z_2^{(n+1,M)}-u_1^{(n+1,L)})\right).
\end{equation*}

Let us recall that now we are considering functionals as in Proposition \ref{finitVese} with $\varphi(s)=s$, $T=I$, and $\Omega=\Omega_i$, $i=1,2$. 
From Proposition \ref{finitVese} and Remark \ref{remarkk} we get that $\xi_1^{(n+1,L)}$, and consequently $u_1^{(n+1,L)}$ is optimal, i.e., {\nnew $-2 \eta_1^{(n+1,L)}\in\partial \hat F_1(\xi_1^{(n+1,L)})$}, if and only if there exists an $M_1^{(n+1)}=(M_{0,1}^{(n+1)},\bar{M}_{1}^{(n+1)})\in V_1\times V_1^d$ with $|\bar{M}_{1}^{(n+1)}(x)|\leq2\alpha$ for all $x\in\Omega_1$ such that
\begin{align}
\langle \bar{M}_{1}^{(n+1)}(x),(\nabla(u_1^{(n+1,L)} + \tilde u_2^{(n)}))(x) \rangle_{\R^d} + 2\alpha \varphi(|(\nabla(u_1^{(n+1,L)} + \tilde u_2^{(n)}))(x)|) = 0 \label{cond1a}\\
-2(u_1^{(n+1,L)}(x)-z_1^{(n+1,L)}(x)) - \operatorname{div} \bar{M}_{1}^{(n+1)}(x) -2 \eta_1^{(n+1,L)}(x) = 0.\label{cond1b} 
\end{align}
for all $x\in\Omega_1$.
Analogously we get that $\xi_2^{(n+1,M)}$, and consequently $u_2^{(n+1,M)}$ is optimal, i.e., {\nnew  $-2 \eta_2^{(n+1,M)}\in\partial \hat F_2(\xi_2^{(n+1,M)})$}, if and only if there exists an $M_2^{(n+1)}=(M_{0,2}^{(n+1)},\bar{M}_{2}^{(n+1)})\in V_2\times V_2^d$ with $|\bar{M}_{2}^{(n+1)}(x)|\leq 2\alpha$ for all $x\in\Omega_2$ such that
\begin{align}
\langle \bar{M}_{2}^{(n+1)}(x),(\nabla(u_1^{(n+1,L)} + u_2^{(n+1,M)}))(x) \rangle_{\R^d} + 2\alpha \varphi(|(\nabla(u_1^{(n+1,L)} + \tilde u_2^{(n+1,M)}))(x)|) &= 0 \label{cond2a}\\
-2(u_2^{(n+1,M)}(x)-z_2^{(n+1,M)}(x)) - \operatorname{div} \bar{M}_{2}^{(n+1)}(x) -2 \eta_2^{(n+1,M)}(x) &= 0, \label{cond2b}
\end{align}
for all $x\in\Omega_2$.
Since $(\bar{M}_1^{(n)}(x))_{n\in\N}$ is bounded for all $x\in\Omega_1$ and $(\bar{M}_2^{(n)}(x))_{n\in\N}$ is bounded for all $x\in\Omega_2$, there exist convergent subsequences $(\bar{M}_1^{(n_k)}(x))_{k\in\N}$ and $(\bar{M}_2^{(n_k)}(x))_{k\in\N}$. Let us denote $\bar{M}_1^{(\infty)}(x)$ and $\bar{M}_2^{(\infty)}(x)$ the {\nnew respective limits} of the sequences. For simplicity we rename such sequences by $(\bar{M}_1^{(n)}(x))_{n\in\N}$ and $(\bar{M}_2^{(n)}(x))_{n\in\N}$.

Note that, by Lemma \ref{eta_bd} (or simply from \eqref{cond1b} and \eqref{cond2b}) the sequences {\nnew $(\eta_1^{(n,L)})_{n\in\N}$ and $(\eta_2^{(n,M)})_{n\in\N}$ are also bounded. Hence there exist convergent subsequences which we denote, for simplicity, again by $(\eta_1^{(n,L)})_{n\in\N}$ and $(\eta_2^{(n,M)})_{n\in\N}$ with limits $\eta_i^{(\infty)}$, $i=1,2$.} By taking in \eqref{cond1a}-\eqref{cond2b} the limits for $n\to \infty$ we obtain

\begin{align*}
\langle \bar{M}_{1}^{(\infty)}(x),(\nabla(u_1^{(\infty)} + u_2^{(\infty)}))(x) \rangle_{\R^d} + 2\alpha \varphi(|(\nabla(u_1^{(\infty)} + u_2^{(\infty)}))(x)|) = 0 \quad \text {for all } x\in\Omega_1\\
-2(u_1^{(\infty)}(x)-z_1^{(\infty)}(x)) - \operatorname{div} \bar{M}_{1}^{(\infty)}(x) -2 \eta_1^{(\infty)}(x) = 0\quad \text {for all } x\in\Omega_1
\end{align*}
\begin{align*}
\langle \bar{M}_{2}^{(\infty)}(x),(\nabla(u_1^{(\infty)} + u_2^{(\infty)}))(x) \rangle_{\R^d} + 2\alpha \varphi(|(\nabla(u_1^{(\infty)} + u_2^{(\infty)}))(x)|) = 0 \quad \text {for all } x\in\Omega_2\\
-2(u_2^{(\infty)}(x)-z_2^{(\infty)}(x)) - \operatorname{div} \bar{M}_{2}^{(\infty)}(x) -2 \eta_2^{(\infty)}(x) = 0\quad \text {for all } x\in\Omega_2
\end{align*}
Since $\operatorname{supp}\eta_1^{(\infty)} = \Gamma_1$ and $\operatorname{supp}\eta_2^{(\infty)} = \Gamma_2$ we have 
\begin{equation}
\begin{split}\label{condOmega1}
\langle \bar{M}_{1}^{(\infty)}(x),(\nabla(u^{(\infty)})(x) \rangle_{\R^d} + 2\alpha \varphi(|(\nabla(u^{(\infty)})(x)|) &= 0 \quad \text {for all } x\in\Omega_1\\
-2T^*((Tu^{(\infty)})(x)-g^{(\infty)}(x)) - \operatorname{div} \bar{M}_{1}^{(\infty)}(x) &= 0\quad \text {for all } x\in\Omega_1\setminus\Gamma_1
\end{split}
\end{equation}
\begin{equation}
\begin{split}\label{condOmega2}
\langle \bar{M}_{2}^{(\infty)}(x),(\nabla(u^{(\infty)})(x) \rangle_{\R^d} + 2\alpha \varphi(|(\nabla(u^{(\infty)})(x)|) &= 0 \quad \text {for all } x\in\Omega_2\\
-2T^*((Tu^{(\infty)})(x)-g^{(\infty)}(x)) - \operatorname{div} \bar{M}_{2}^{(\infty)}(x) &= 0\quad \text {for all } x\in\Omega_2\setminus\Gamma_2.
\end{split}
\end{equation}

Observe now that from Proposition \ref{finitVese} we also have that $0\in\mathcal{J}(u^{(\infty)})$ if and only if there exists $M^{(\infty)} = (M_0^{(\infty)},\bar{M}^{(\infty)})$ with $|\bar{M}_{0}^{(\infty)}(x)|\leq 2\alpha$ for all  $x\in\Omega$ such that 
\begin{equation}
\begin{split}\label{condOmega}
\langle \bar{M}^{(\infty)}(x),(\nabla(u^{(\infty)})(x) \rangle_{\R^d} + 2\alpha \varphi(|(\nabla(u^{(\infty)})(x)|) &= 0 \quad \text {for all } x\in\Omega\\
-2T^*((Tu^{(\infty)})(x)-g^{(\infty)}(x)) - \operatorname{div} \bar{M}^{(\infty)}(x) &= 0\quad \text {for all } x\in\Omega.
\end{split}
\end{equation}

Note that  $\bar{M}_j^{(\infty)}(x),j=1,2$, for $x\in\Omega_1\cap\Omega_2$ satisfies both \eqref{condOmega1} and \eqref{condOmega2}. Hence let us choose 
$$
M^{(\infty)}(x)=
\begin{cases} 
M_1^{(\infty)}(x) \quad \text{if } x\in\Omega_1\setminus\Gamma_1  \\
M_2^{(\infty)}(x) \quad \text{if } x\in(\Omega_2\setminus\Omega_1)\cup\Gamma_1 
\end{cases}.
$$
With this choice of $M^{(\infty)}$  equations \eqref{condOmega1} - \eqref{condOmega} are valid and hence $u^{(\infty)}$ is optimal in $\Omega$. 
\end{proof}

\begin{remark}
\begin{enumerate}
\item[(i)]
If $\nabla u^{(\infty)}(x)\not=0$ for $x\in\Omega_j$, $j=1,2$, then $\bar{M}_j^{(\infty)}$ is given as in equation \eqref{Mbar} by 
{\nnew $$
\bar{M}_j^{(\infty)}(x)=-2\alpha \frac{(\nabla u^{(\infty)}\mid_{\Omega_j})(x)}{|(\nabla u^{(\infty)}\mid_{\Omega_j})(x)|}.
$$}
{\nnew \item[(ii)] The boundedness of the sequences $(\tilde u_1^{(n)})_{n\in \mathbb{N}}$ and $(\tilde u_2^{(n)})_{n\in \mathbb{N}}$  has been technically used for showing  the existence of an optimal decomposition $u^{(\infty)}=u^{(\infty)}_1+ u^{(\infty)}_2$ in the proof of Theorem \ref{convergence:proof}. Their boundedness is guaranteed as in Lemma \ref{bdtu12} by the use of the partition of the unity $\{\chi_1, \chi_2\}$. Let us emphasize that there is no way of obtaining the boundedness of the local sequences $(u_1^{(n,L)})_{n\in \mathbb{N}}$ and $(u_2^{(n,M)})_{n\in \mathbb{N}}$ otherwise. In Figure \ref{fig:1Dnum} we show that the local sequences can become unbounded in case we do not modify them by means of the partition of the unity.
}
\item[(iii)]Note that for deriving the optimality condition  \eqref{condOmega} for $u^{(\infty)}$ we combined the respective conditions \eqref{condOmega1} and \eqref{condOmega2} for $u_1^{(\infty)}$ and $u_2^{(\infty)}$. In doing that, we {\nnew strongly took} advantage of the overlapping property of the subdomains, hence avoiding a fine analysis of $\eta_1^{(\infty)}$ and $\eta_2^{(\infty)}$ on the interfaces $\Gamma_1$ and $\Gamma_2$. This is the major advantage of this analysis with respect to the one provided in \cite{FS} for nonoverlapping domain decompositions.
\end{enumerate}
\end{remark}

\section{A parallel algorithm and its convergence}\label{Convsec2}
The parallel version of the previous algorithm \eqref{schw_sp:it2} reads as follows: Pick an initial $V_1 + V_2 \ni \tilde u_1^{(0)}+ \tilde u_2^{(0)} : = u^{(0)} \in \mathcal{H}$, for example $\tilde u_i^{(0)}=0$, $i=1,2$, and iterate
\begin{equation}
\label{parallel_schw_sp:it2}
\left \{ 
\begin{array}{ll}
\left \{ 
\begin{array}{ll}
u_1^{(n+1,0)} = \tilde u_1^{(n)}&\\
u_1^{(n+1,\ell+1)} =  \operatorname{argmin}_{\stackrel{u_1 \in V_1}{\operatorname{Tr}\mid_{\Gamma_1}u_1=0}} \mathcal J_1^s(u_1+ \tilde u_2^{(n)}, u_1^{(n+1,\ell)}) & \ell=0,\dots, L-1\\
\end{array}\right. &\\
\left \{ 
\begin{array}{ll}
u_2^{(n+1,0)} = \tilde{u}_2^{(n)}&\\
u_2^{(n+1,m+1)} = \operatorname{argmin}_{\stackrel{u_2 \in V_2}{\operatorname{Tr}\mid_{\Gamma_2}u_2=0}} \mathcal J_2^s(\tilde u_1^{(n)}+ u_2, u_2^{(n+1,m)}) &m=0,\dots, M -1\\
\end{array}\right. &\\
u^{(n+1)}:=\frac{u_{1}^{(n+1,L)} + u_{2}^{(n+1,M)} + u^{(n)}}{2}\\
\tilde{u}_1^{(n+1)}:=\chi_1\cdot u^{(n+1)} \\
\tilde{u}_2^{(n+1)}:=\chi_2\cdot u^{(n+1)}
\end{array}
\right .
\end{equation} 

We are going to propose similar convergence results as for the sequential algorithm.

\begin{proposition}[Convergence properties]
\label{parallel_weak-conv}
Let us assume that $\|T\|<1$. The parallel algorithm (\ref{parallel_schw_sp:it2}) produces a sequence $(u^{(n)})_{n\in \mathbb{N}}$ in  $\mathcal{H}$ with the following properties:
\begin{itemize}
\item[(i)] $\mathcal{J}(u^{(n)}) > \mathcal{J}(u^{(n+1)})$ for all $n \in \mathbb{N}$ (unless $u^{(n)}= u^{(n+1)}$);
\item[(ii)] $\lim_{n \to \infty} \| u^{(n+1)} -  u^{(n)}\|_{2} = 0$;
\item[(iii)] the sequence  $(u^{(n)})_{n\in \mathbb{N}}$ has subsequences which converge in $\mathcal{H}$.
\end{itemize}
\end{proposition}

\begin{proof}
With the same argument as in the proof of Theorem \ref{weak-conv}, we obtain
$$
\mathcal {J}(u^{(n)}) -  \mathcal {J}(u_1^{(n+1,L)} + \tilde u_2^{(n)})\geq C \sum_{\ell=0}^{L-1}\| u_1^{(n+1,\ell+1)} - u_1^{(n+1,\ell)}\|_{2}^2
$$
and
$$
 \mathcal {J}(u^{(n)}) -  \mathcal {J}(\tilde u_1^{(n)}   + u_2^{(n+1,M)})\geq C  \sum_{m=0}^{M-1}\| u_2^{(n+1,m+1)} - u_2^{(n+1,m)}\|_{2}^2.
$$
Hence, by summing and halving  
\begin{eqnarray*}
&& \mathcal {J}(u^{(n)}) -  \frac{1}{2}( \mathcal {J}(u_1^{(n+1,L)}   + \tilde u_2^{(n)}) + \mathcal {J}(\tilde u_1^{(n)}   + u_2^{(n+1,M)})) \\
&\geq& \frac{C}{2} \left ( \sum_{\ell=0}^{L-1}\| u_1^{(n+1,\ell+1)} - u_1^{(n+1,\ell)}\|_{2}^2 + \sum_{m=0}^{M-1}\| u_2^{(n+1,m+1)} - u_2^{(n+1,m)}\|_{2}^2 \right ).
\end{eqnarray*}
We recall that $\mathcal{J}(u^{(n)})=\|Tu^{(n)}-g\|_2^2 + 2\alpha|\nabla u^{(n)}|(\Omega)$ (and $T$ is linear). {\nnew Then, by the standard inequality $(a^2+b^2) \geq \frac{1}{2} (a+b)^2$ for $a,b>0$, 
we have}
\begin{eqnarray*}
\left \| T u^{(n+1)}- g \right\|^2_{2} &=& \left \| T \left (\frac{(u_1^{(n+1,L)} + u_2^{(n+1,M)}) + u^{(n)}}{2} \right ) - g \right\|^2_{2} \\
&\leq & \frac{1}{2} \| T(u^{(n+1,L)}_1 +  \tilde u^{(n)}_2 ) -g \|^2_{2} + \frac{1}{2} \| T(\tilde u^{(n)}_1 + u^{(n+1,M)}_2 ) -g \|^2_{2}.
\end{eqnarray*}
Moreover we have
\begin{eqnarray*}
|\nabla(u^{(n+1)})|(\Omega) &\leq& \frac{1}{2} \left ( |\nabla(u_1^{(n+1,L)} + \tilde u_2^{(n)})|(\Omega) + |\nabla(\tilde u^{(n)}_1 + u_2^{(n+1,M)})|(\Omega) \right ).
\end{eqnarray*}
By the last two inequalities we immediately show that
$$
\mathcal  J (u^{(n+1)}) \leq \frac{1}{2} \left (\mathcal {J}(u_1^{(n+1,L)} + \tilde u_2^{(n)})+ \mathcal {J}(\tilde u^{(n)}_1   + u_2^{(n+1,M)}) \right ),
$$
hence
$$
\mathcal {J}(u^{(n)})- \mathcal  J (u^{(n+1)})$$\begin{equation}
\label{coerc3}\geq \frac{C}{2} \left( \sum_{\ell=0}^{L-1}\| u_1^{(n+1,\ell+1)} - u_1^{(n+1,\ell)}\|_{2}^2 +\sum_{m=0}^{M-1}\| u_2^{(n+1,m+1)} - u_2^{(n+1,m)}\|_{2}^2 \right) \geq 0.
\end{equation}
Since the sequence $(\mathcal {J}(u^{(n)}))_{n \in \mathbb{N}}$ is monotonically decreasing and bounded from below by 0, it is also convergent. From (\ref{coerc3}) and the latter convergence we deduce 
\begin{equation}
\label{asymp_reg3}
\left ( \sum_{\ell=0}^{L-1}\| u_1^{(n+1,\ell+1)} - u_1^{(n+1,\ell)}\|_{\mathcal H}^2 + \sum_{m=0}^{M-1}\| u_2^{(n+1,m+1)} - u_2^{(n+1,m)}\|_{2}^2 \right ) \rightarrow 0, \quad n \to \infty.
\end{equation}
In particular, {\nnew by again using} $(a^2+b^2) \geq \frac{1}{2} (a+b)^2$ for $a,b>0$ and the triangle inequality, {\nnew we also have
\begin{equation}
\label{asymp_reg4}
\| u^{(n)} -  u^{(n+1)} \|_{2} \rightarrow 0, \quad n\to \infty.
\end{equation}}
The rest of the proof follows analogous arguments as in that of Proposition \ref{weak-conv}.
\end{proof}

{\nnew Analogous results as the one stated in Lemma \ref{bdtu12} and Lemma \ref{eta_bd} also hold in the parallel case.}
%
With these preliminary results the following theorem holds:

\begin{theorem}[Convergence to minimizers]
Assume $\|T\|<1$. Then accumulation points of the sequence $(u^{(n)})_{n\in \mathbb{N}}$ produced by algorithm (\ref{parallel_schw_sp:it2}) are minimizers of $\mathcal{J}$. If $\mathcal{J}$ has a unique minimizer then the sequence $(u^{(n)})_{n\in \mathbb{N}}$ converges to it.
\end{theorem}

\begin{proof}
Note that $u^{(n+1)}$ is the average of the current iteration and the previous, i.e.,
$$
u^{(n+1)} = \frac{u_1^{(n+1,L)} + u_2^{(n+1,M)} + u^{(n)}}{2}.
$$
Observe that the sequences $(u_1^{(n+1,L)})_{n\in\mathbb{N}}$, $(u_2^{(n+1,M)})_{n\in\mathbb{N}}$ and $(u^{(n)})_{n\in\mathbb{N}}$ are bounded. Hence there exist convergent subsequences. By taking the limit for $n\to\infty$ we obtain
$$
u^{(\infty)} = \frac{u_1^{(\infty)} + u_2^{(\infty)} + u^{(\infty)}}{2}
$$
which is equivalent to 
$$
u^{(\infty)} = u_1^{(\infty)} + u_2^{(\infty)}.
$$
With this observation the rest of the proof follows analogous arguments as in that of Theorem \ref{convergence:proof}.
\end{proof}


\section{Applications and Numerics}\label{Numersec}

 In this section we shall present the application of the sequential algorithm \eqref{schw_sp} for the minimization of $\mathcal J$ in one and two dimensions. In particular, we show how to implement the dual method of Chambolle \cite{Ch} in order to compute the orthogonal projection $P_{\alpha K}(g)$ in the oblique thresholding, and we give a detailed explanation of the domain decompositions used in the numerics. Furthermore we present numerical examples for image {\it inpainting}, i.e., the recovery of missing parts of images by minimal total variation interpolation, and compressed sensing \cite{cand,carotaXX,cataXX,do04}, the nonadaptive compressed acquisition of images for a classical toy problem inspired by magnetic resonance imaging (MRI) \cite{carotaXX,ludopa07}. The numerical examples of this section and respective Matlab codes can be found at \cite{WEB}.

\subsection{Computation of $P_{\alpha K}(g)$}
To solve the subiterations in \eqref{schw_sp} we compute the minimizer by means of oblique thresholding. More precisely, let us denote $u_2=\tilde{u}_2^{(n)}$, $u_1=u_1^{(n+1,\ell+1)}$, and $z_1=u_1^{(n+1,\ell)}+ T^*( g - T u_2 - T u_1^{(n+1,\ell)})$. We shall compute the minimizer $u_1$ of the first subminimization problem by
$$
u_1 = (I- P_{\alpha K})(z_1+ u_2 -\eta)-u_2 \in V_1
$$
for an $\eta \in V_1$ with $\operatorname{supp}\eta = \Gamma_1$ which fulfills
$$
\operatorname{Tr}\mid_{\Gamma_1} (\eta) = \operatorname{Tr}\mid_{\Gamma_1} \left(z_1 + P_{\alpha K} (\eta - z_1 - u_2)\right).
$$
Hence the element $\eta\in V_1$ is a limit of the corresponding fixed point iteration
\begin{equation}
\eta^{(0)}\in V_1,\supp \eta^{(0)}=\Gamma_1, \quad  \eta^{(m+1)} = (\operatorname{Tr}\mid_{\Gamma_1})^*\operatorname{Tr}\mid_{\Gamma_1} \left(z_1 + P_{\alpha K} (\eta^{(m)} - z_1 - u_2)\right),\; m\geq 0.
\end{equation}
Here $K$ is defined as in Section \ref{notations}, i.e., 
$$
K=\left\{\operatorname{div } p: p\in \mathcal{H}^d, \left|p(x)\right|_\infty\leq 1 \quad\forall x \in\Omega \right\}.
$$
To compute the projection onto $\alpha K$ in the oblique thresholding we use an algorithm proposed by Chambolle in \cite{Ch}. { His algorithm is based on considerations of the convex conjugate of the total variation and on exploiting the corresponding optimality condition. It amounts to compute $P_{\alpha K}(g)$ approximately by $\alpha\operatorname{div} p^{(n)}$, where $p^{(n)}$ is the $nth$ iterate} of the following semi-implicit gradient descent algorithm:

\begin{quote}
Choose $\tau>0$, let $p^{(0)}=0$ and, for any $n\geq 0$, iterate
\begin{eqnarray*}
p^{(n+1)}(x) = \frac{p^{(n)}(x) + \tau (\nabla(\operatorname{div} p^{(n)}-g/\alpha))(x)}{1+\tau \left|(\nabla(\operatorname{div} p^{(n)}-g/\alpha))(x)\right|}.
\end{eqnarray*}
\end{quote}

For $\tau>0$ sufficiently small, i.e., $\tau< 1/8$, the iteration $\alpha\operatorname{div} p^{(n)}$ was shown to converge to $P_{\alpha K}(g)$ as $n\rightarrow\infty$ (compare \cite[Theorem 3.1]{Ch}). {\nnew Let us stress that we propose here this algorithm just for the ease of its presentation; its choice for the approximation of projections is of course by no means a restriction and one may want to implement other recent, and perhaps faster strategies, e.g., \cite{ChDa,DaSig,GoOs,breg,WBA}.}

\subsection{Domain decompositions}
\label{tricks}

\begin{figure}[htbp]
\begin{center}
\includegraphics[height=5cm]{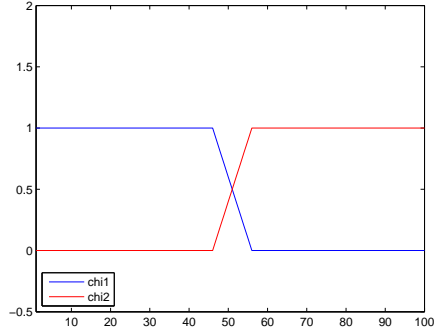}
\end{center}    
\caption{\small Auxiliary functions $\chi_1$ and $\chi_2$ for an overlapping domain decomposition with two subdomains.}
\label{chi}
\end{figure}

\begin{figure}[htbp]

\begin{picture}(-140,20)(-208,10)

\put(35,23){$\Omega_2$}
\put(-22,+7){\textcolor{red}{$\Gamma_2$}}
\put(-20,20){\line(1,0){120}}
\put(-20,22){\line(0,-1){4}}
\put(100,22){\line(0,-1){4}}

\put(-20,0){\textcolor{red}{\circle{5}}}
\put(20,0){\textcolor{red}{\circle{5}}}
\put(-100,0){\textcolor{blue}{\thicklines\line(1,0){200}}}
\put(-100,-2){\textcolor{blue}{\thicklines\line(0,1){4}}}
\put(100,-2){\textcolor{blue}{\thicklines\line(0,1){4}}}

\put(18,-15){\textcolor{red}{$\Gamma_1$}}
\put(-100,-20){\line(1,0){120}}
\put(20,-22){\line(0,1){4}}
\put(-100,-22){\line(0,1){4}}
\put(-45,-29){$\Omega_1$}
\end{picture}
$\vspace{1cm}$
\caption{\small Overlapping domain decomposition in 1D.}
\label{fig1Ddd}
\end{figure}
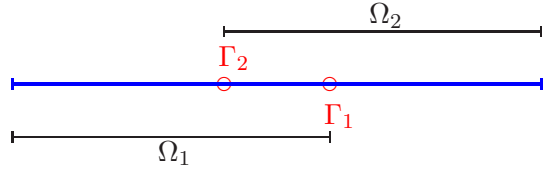

In one dimension the domain $\Omega=[a,b]$ is split into two overlapping intervals. Let $\left|\Omega_1\cap\Omega_2\right|=:G$ be the size of the overlap of $\Omega_1$ and $\Omega_2$. Then we set $\left|\Omega_1\right|=: n_1 = \left\lceil \frac{N+G}{2}\right\rceil$,  $\Omega_1=[a,n_1]$ and $\Omega_2=[n_1-G+1,b]$. The interfaces $\Gamma_1$ and $\Gamma_2$ are located in $i=n_1+1$ and $n_1-G$ respectively (cf. Figure \ref{fig1Ddd}). The auxiliary functions $\chi_1$ and $\chi_2$ can be chosen in the following way (cf. Figure \ref{chi}):
\begin{align*}
&\chi_1(x_i)=
\begin{cases}
1 \quad &x_i \in \Omega_1\setminus\Omega_2\\
1-\frac{1}{G}(i-(n_1-G+1)) \quad &x_i\in\Omega_1\cap\Omega_2
\end{cases}\\
&\chi_2(x_i)=
\begin{cases}
1 \quad &x_i \in \Omega_2\setminus\Omega_1\\
\frac{1}{G}(i-(n_1-G+1)) \quad &x_i\in\Omega_1\cap\Omega_2
\end{cases}.
\end{align*}
Note that $\chi_1(x_i) + \chi_2(x_i)=1$ for all $x_i\in\Omega$ (i.e for all $i=1,\ldots,N$).

\begin{center}
\begin{figure}[htbp]
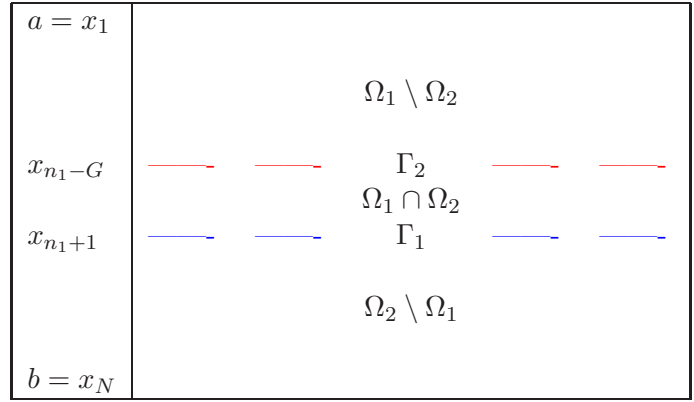
\begin{center}
\begin{tabular}{|l| p{1 cm} p{1 cm} c p{1 cm} p{1 cm} |}
\hline
$a=x_1$ & & & & &\\
 & & & & &\\
 & & & $\Omega_1\setminus\Omega_2$ & &\\
 & & & & &\\
$x_{n_1-G}$&\textcolor{red}{ ------- }&\textcolor{red}{ ------- }& $\Gamma_2$ &\textcolor{red}{ ------- }&\textcolor{red}{ ------- }\\
 & & & $\Omega_1\cap\Omega_2$ & &\\
$x_{n_1+1}$&\textcolor{blue}{ ------- }&\textcolor{blue}{ ------- }& $\Gamma_1$ &\textcolor{blue}{ ------- }&\textcolor{blue}{ ------- }\\
 & & & & &\\
 & & & $\Omega_2\setminus\Omega_1$ & &\\
 & & & & &\\
$b=x_N$ & & & & &\\
\hline
\end{tabular}\end{center}
\caption{Decomposition of the image in two domains $\Omega_1$ and $\Omega_2$.}
\label{figill}
\end{figure}
\end{center}

In two dimensions the domain $\Omega=[a,b]\times[c,d]$ is split in an analogous way with respect to its rows. In particular we have $\Omega_1=[a,n_1]\times[c,d]$ and $\Omega_2=[n_1-G+1,b]\times[c,d]$, compare Figure \ref{figill}. The splitting in more than two domains is done similarly: 
\begin{quote}
Set $\Omega=\Omega_1\cup\ldots\cup\Omega_\mathcal N$, the domain $\Omega$ decomposed into $\mathcal N$ domains $\Omega_i$, $i=1,\ldots,\mathcal N$, where $\Omega_i$ and $\Omega_{i+1}$ are overlapping for $i=1,\ldots, \mathcal N-1$. Let $|\Omega_i\cap\Omega_{i+1}|=:G$ equidistant for every $i=1,\ldots,\mathcal N-1$. Set $s=\left\lceil N_1/\mathcal N\right\rceil$. Then
\begin{eqnarray*}
& & \Omega_1=[1,s+\frac{G}{2}]\times[c,d]\\
& & \textrm{for } i=2:\mathcal N-1\\
& & \quad\Omega_i=[(i-1)s-\frac{G}{2}+1,is + \frac{G}{2}]\times[c,d]\\
& & \textrm{end}\\
& & \Omega_\mathcal N=[(\mathcal N-1)s-\frac{G}{2}+1,N_1]\times[c,d].
\end{eqnarray*}
\end{quote}
The auxiliary functions $\chi_i$ can be chosen in an analogous way as in the one dimensional case:
$$
\chi_i(x_{i_1},y_{i_2})=
\begin{cases}
\frac{1}{G}(i_1-((i-1)s-G/2+1)) \quad &(x_{i_1},y_{i_2})\in\Omega_{i-1}\cap\Omega_i\\
1 \quad &(x_{i_1},y_{i_2}) \in \Omega_i\setminus(\Omega_{i-1}\cup\Omega_{i+1})\\
1-\frac{1}{G}(i_1-(is-G/2+1)) \quad &(x_{i_1},y_{i_2})\in\Omega_{i}\cap\Omega_{i+1}
\end{cases}
$$
for $i=1,\ldots,\mathcal{N}$ with $\Omega_0=\Omega_{\mathcal{N}+1}=\emptyset$.


To compute the fixed point $\eta$ of \eqref{fixpt} in an efficient way we make the following considerations, which allow to restrict the computation from $\Omega_1$ to a relatively small stripe around the interface. The fixed point $\eta$ is actually supported on $\Gamma_1$ only, i.e., $\eta(x)=0$ in $\Omega_1 \setminus \Gamma_1$. Hence, we restrict the fixed point iteration for $\eta$ to a relatively small stripe $\hat{\Omega}_1\subset\Omega_1$
Analogously, one implements the minimizations of $\eta_2$ on $\hat{\Omega}_2$. A similar trick was also used in  \cite{FS} to compute suitable Lagrange multipliers at the interfaces of the nonoverlapping domains. However, there we needed to consider larger ``bilateral stripes'' around the support of the multiplier, making the numerical computation slightly more demanding for that algorithm.

\subsection{Numerical experiments}
 In the following we present numerical examples for the sequential algorithm \eqref{schw_sp:it2} in two particular applications: signal interpolation/image inpainting, and { compressed sensing}.

\begin{figure}[htbp]
\begin{center}
    \subfigure[]{\label{l1}\includegraphics[width=7.5cm]{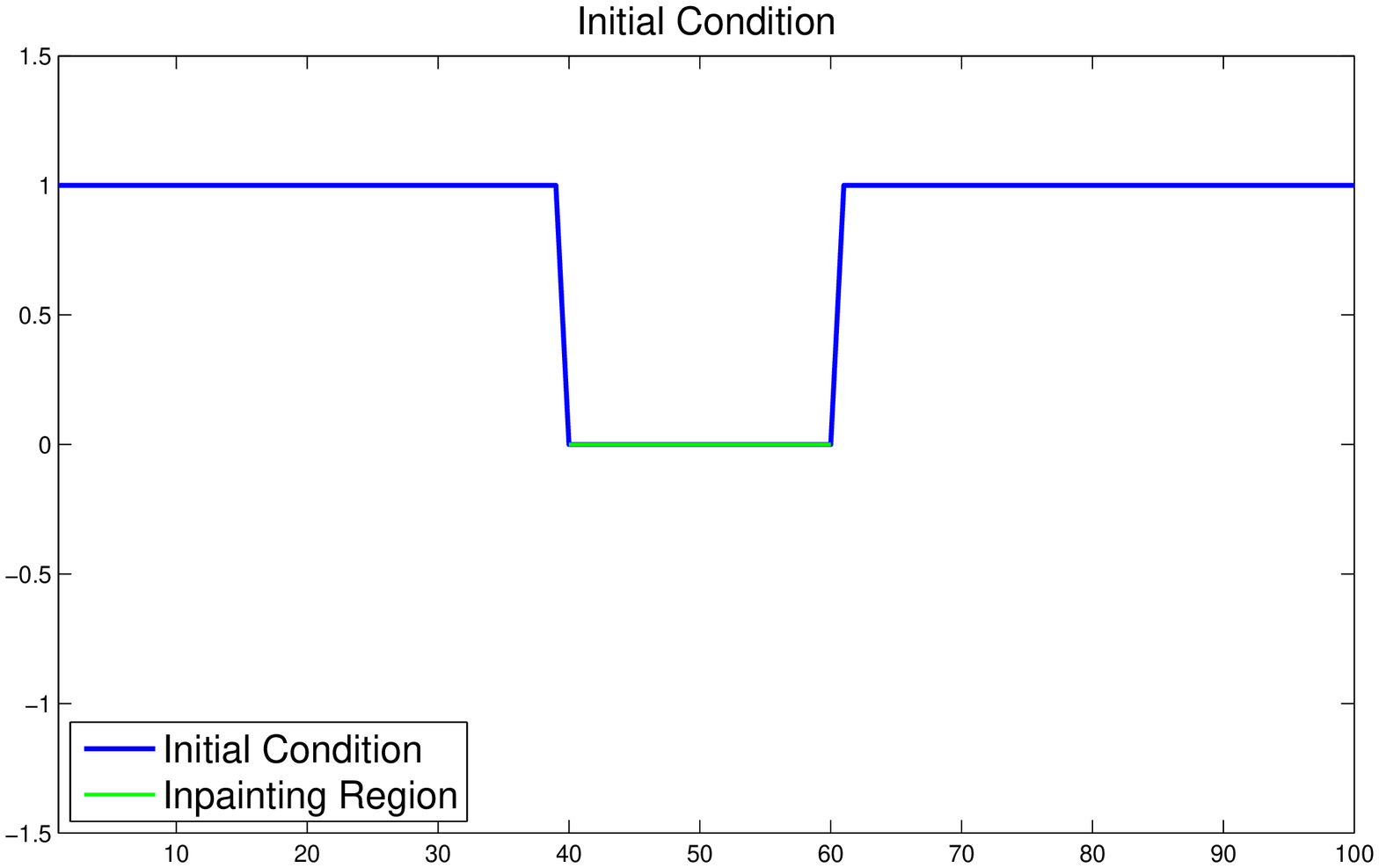}}
    \subfigure[]{\label{l2}\includegraphics[width=7.5cm]{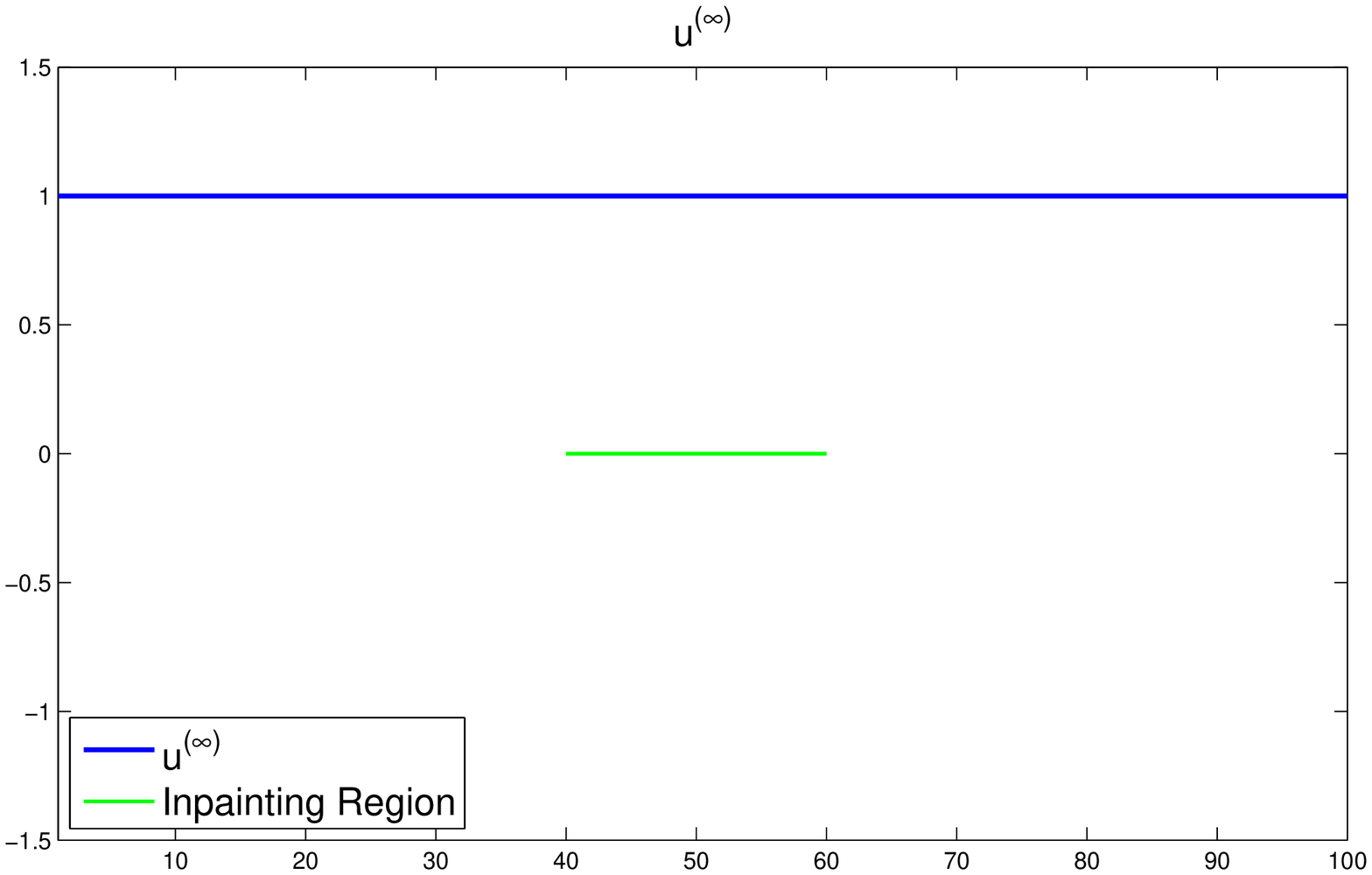}}
\\
 \subfigure[]{\label{l2}\includegraphics[width=7.5cm]{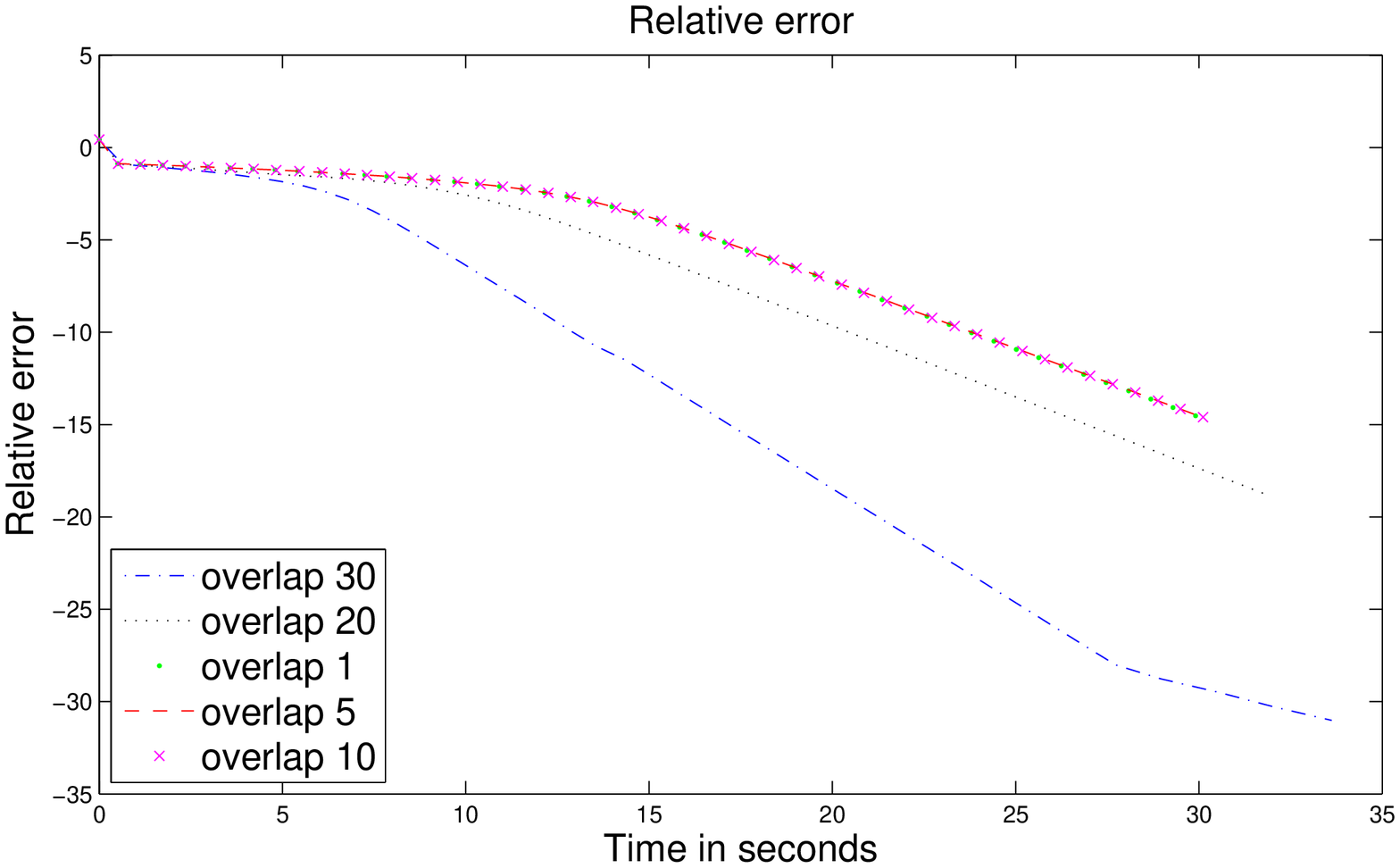}}
 \subfigure[]{\label{l2}\includegraphics[width=7.5cm]{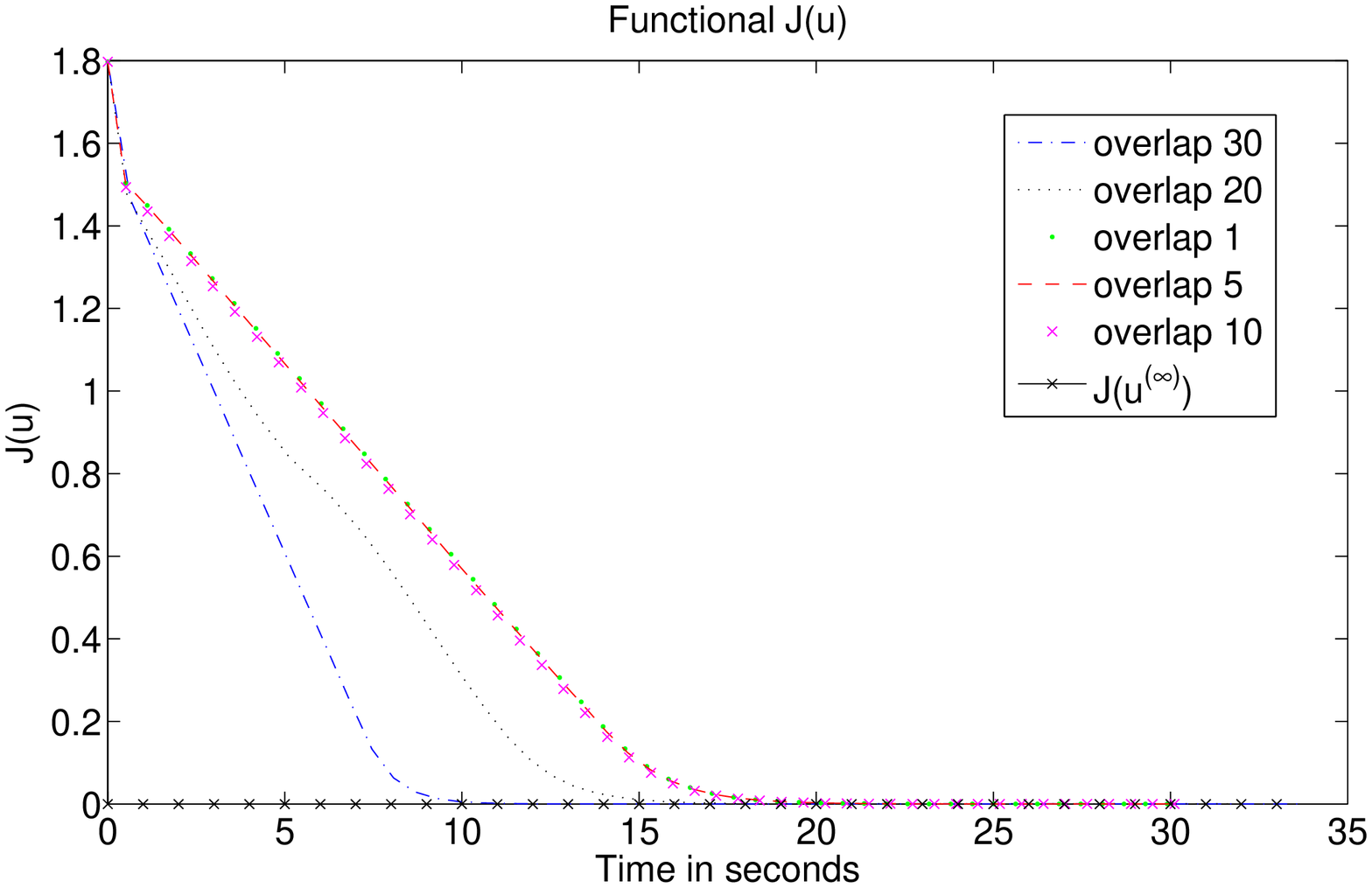}}
\end{center}    
\caption{\small \nnew We present a numerical experiment related to the interpolation of a 1D signal by total variation minimization. The original signal is only provided outside of the green subinterval. The initial datum $g$ is shown in (a). As expected, the minimizer $u^{(\infty)}$ is the constant vector $1$, as shown in (b). In (c) and (d) we display the rates of decay of the relative error and of the value of $\mathcal J$ respectively, for applications of the algorithm \eqref{schw_sp:it2} with different sizes G=1,5,10,20,30 of the overlapping region of two subintervals.}
\label{fig:1Dnum1}
\end{figure}

\begin{figure}[htbp]
\begin{center}
    \subfigure[]{\label{l1}\includegraphics[width=7.5cm]{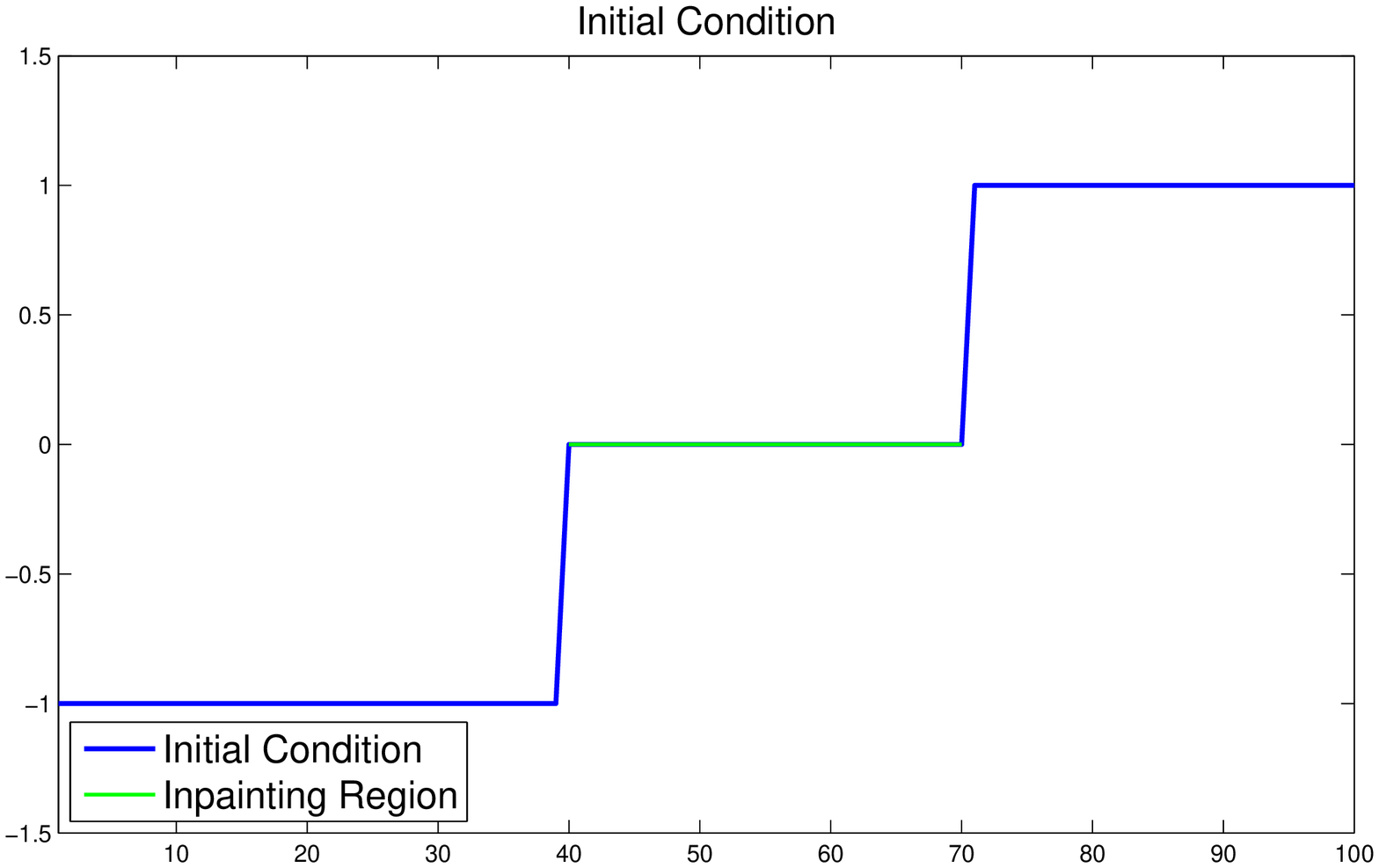}}
    \subfigure[]{\label{l2}\includegraphics[width=7.5cm]{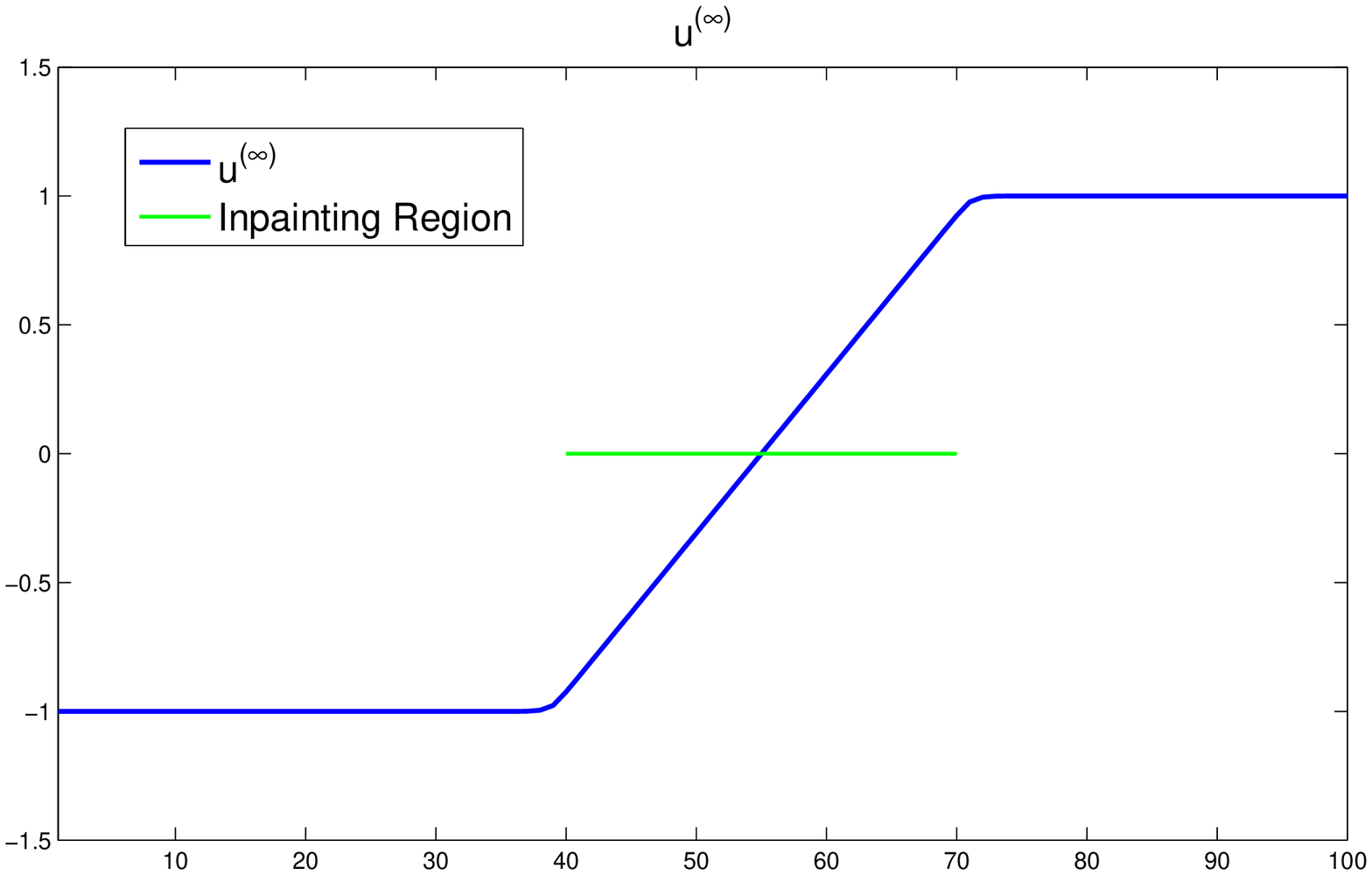}}
\\
 \subfigure[]{\label{l2}\includegraphics[width=7.5cm]{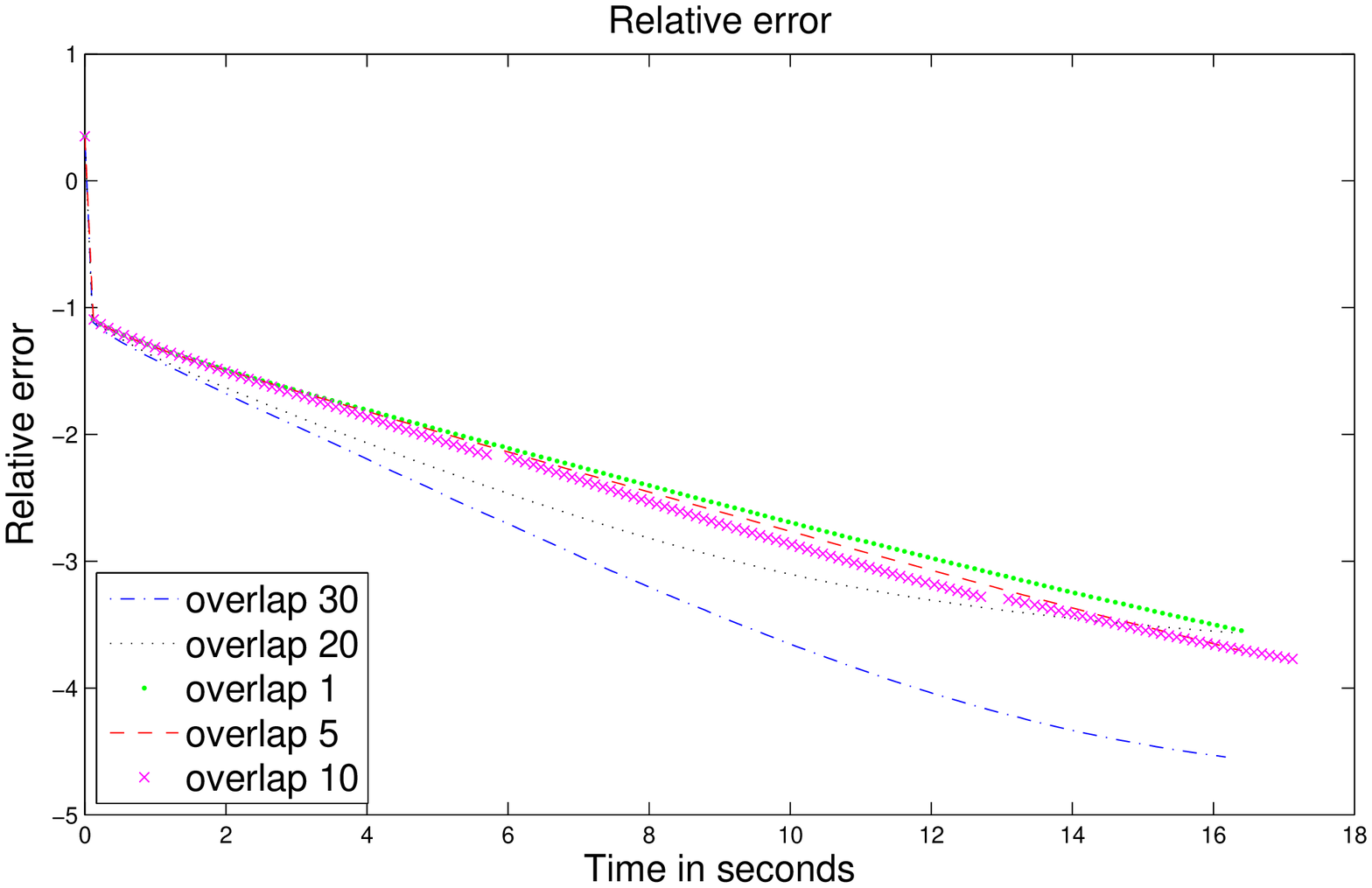}}
 \subfigure[]{\label{l2}\includegraphics[width=7.5cm]{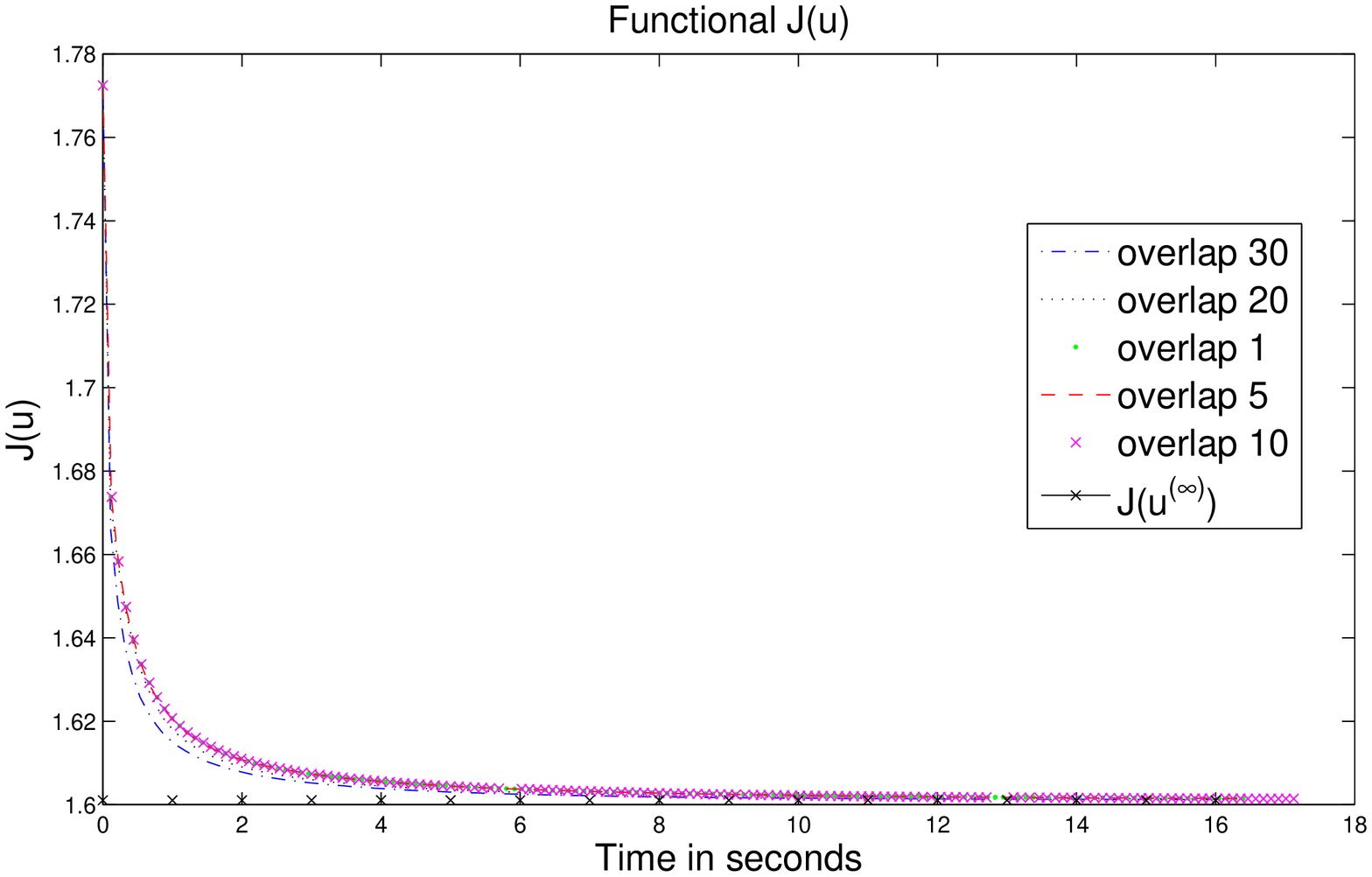}}
\end{center}    
\caption{\small \nnew We show a second example of total variation interpolation in 1D. The initial datum $g$ is shown in (a). As expected, a minimizer $u^{(\infty)}$ is (nearly) a piecewise linear function, as shown in (b). In (c) and (d) we display the rates of decay of the relative error and of the value of $\mathcal J$ respectively, for applications of the algorithm \eqref{schw_sp:it2} with different sizes G=1,5,10,20,30 of the overlapping region of two subintervals.
}
\label{fig:1Dnum2}
\end{figure}

\begin{figure}[htbp]
\begin{center}
    \subfigure[]{\label{l1}\includegraphics[height=5.5cm]{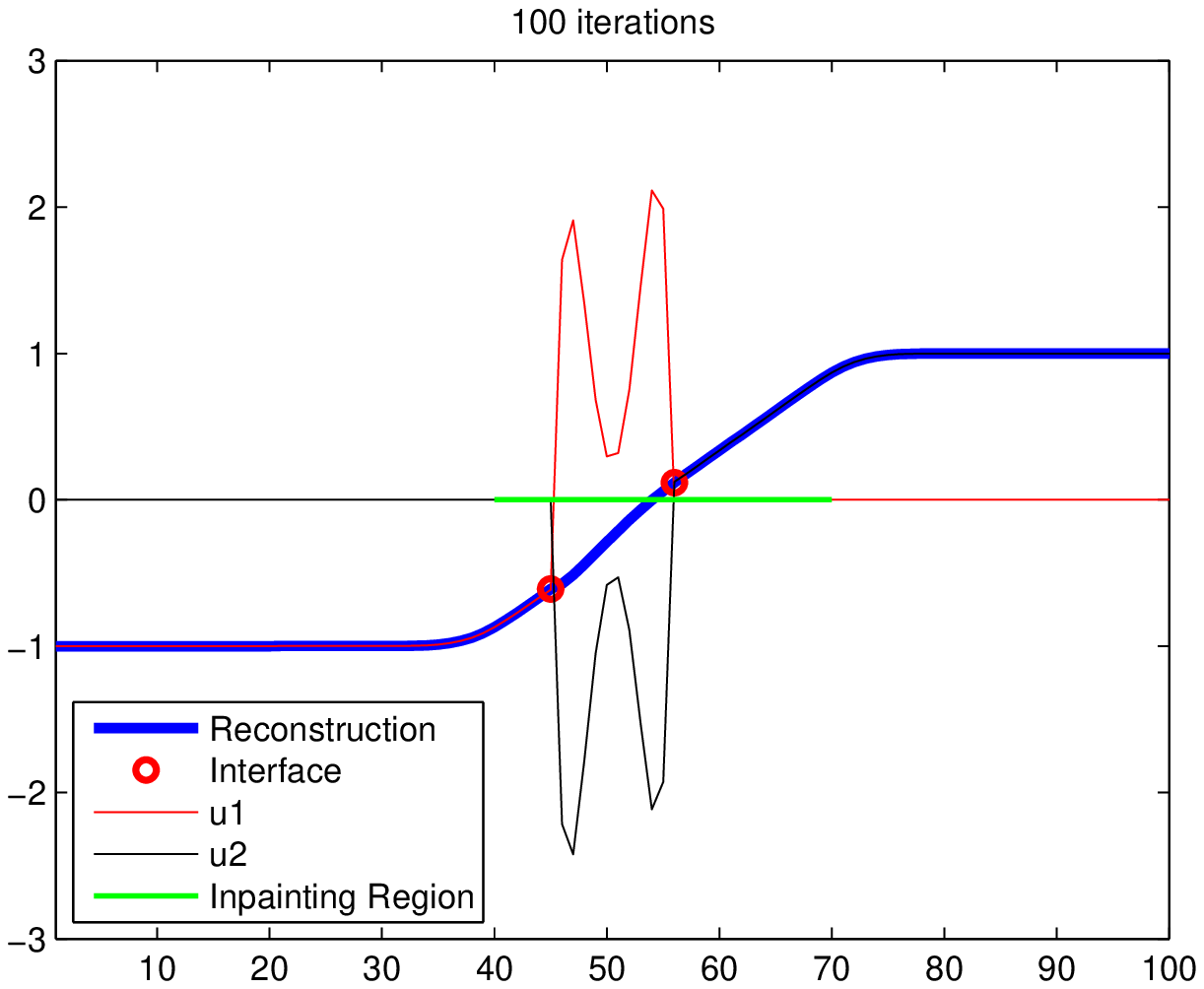}}
    \subfigure[]{\label{l2}\includegraphics[height=5.5cm]{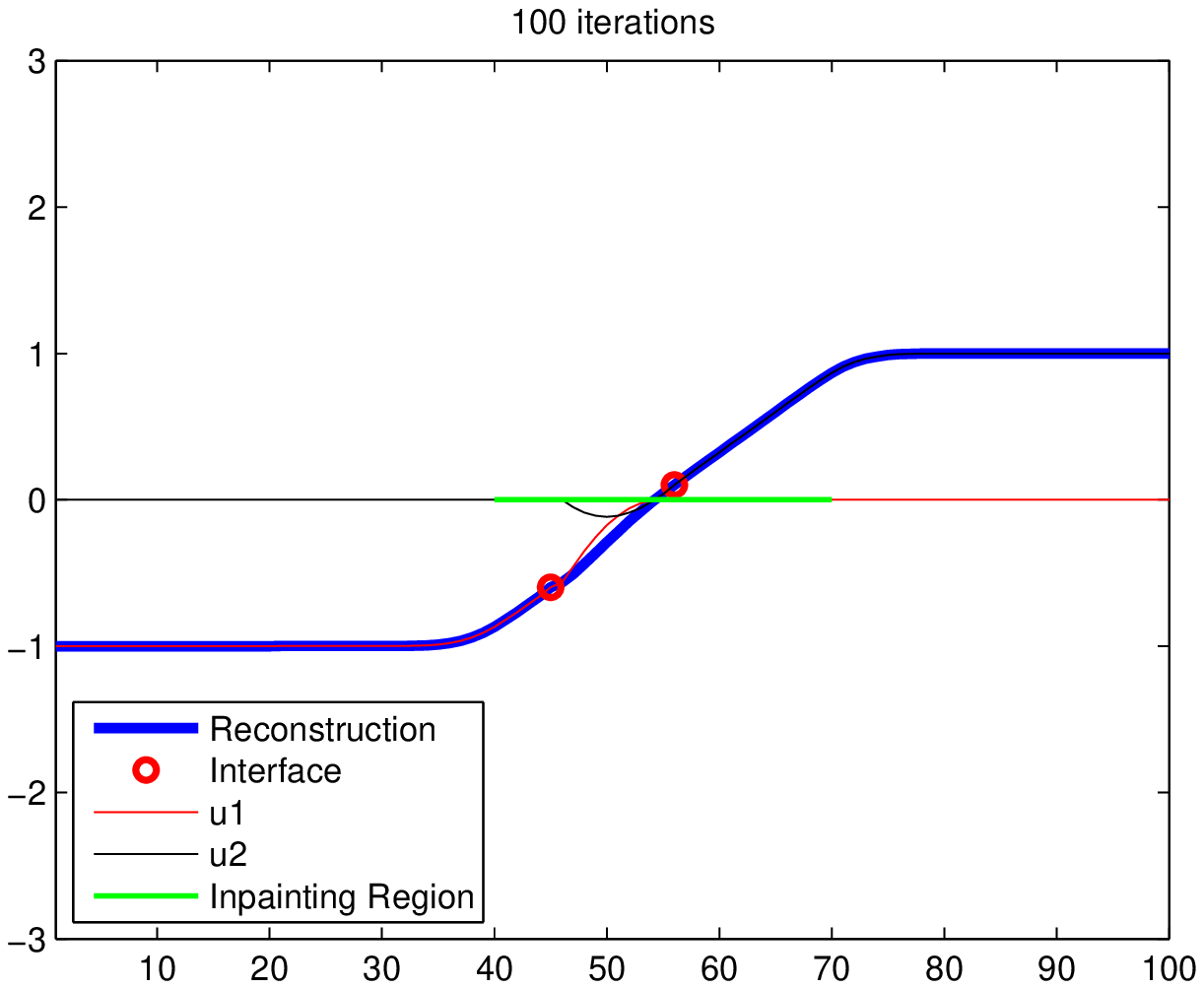}}
\end{center}    
\caption{\small Here we present two numerical experiments related to the interpolation of a 1D signal by total variation minimization. The original signal is only provided outside of the green subinterval. On the left we show an application of algorithm \eqref{schw_sp:it2} when no correction with the partition of unity is provided. In this case, the sequence of the local iterations $u_1^{(n)}, u_2^{(n)}$ is unbounded. On the right we show an application of algorithm \eqref{schw_sp:it2} with the use of the partition of unity which enforces the uniform boundedness of the local iterations $u_1^{(n)}, u_2^{(n)}$.}
\label{fig:1Dnum}
\end{figure}

{\nnew In Figure \ref{fig:1Dnum1} and  Figure \ref{fig:1Dnum2} we show a partially corrupted 1D signal on an interval $\Omega$ of $100$ sampling points, with a loss of information on an interval $D \subset \Omega$. The domain $D$ of the missing signal points is marked with green. These signal points are reconstructed by total variation interpolation, i.e., minimizing the functional $\mathcal J$ in \eqref{functotalfin} with $\alpha=0.4$ and $Tu  = 1_{\Omega\setminus D} \cdot u$, where $1_{\Omega\setminus D}$ is the indicator function of $\Omega\setminus D$. A minimizer $u^{(\infty)}$ of $\mathcal J$ is precomputed with an algorithm working on the whole interval $\Omega$ without any decomposition. We show also the decay of relative error and of the value of the energy $\mathcal J$ for applications of algorithm \eqref{schw_sp:it2} on two subdomains and with different overlap sizes $G=1,5,10,20,30$. The fixed points $\eta$'s are computed on a small interval $\hat{\Omega}_i$, $i=1,2$, of size $2$.
These results confirm the behavior of the algorithm \eqref{schw_sp:it2} as predicted by the theory; the algorithm monotonically decreases $\mathcal J$ and computes a minimizer, independently of the size of the overlapping region. A larger overlapping region does not necessarily imply a slower convergence. In these figures we do compare the speed in terms of CPU time.
In Figure \ref{fig:1Dnum} we also illustrate the effect of implementing the BUPU within the domain decomposition algorithm. In this case, with datum $g$ as in Figure \ref{fig:1Dnum2}, we chose $\alpha=1$ and an overlap of size $G=10$. The fixed points $\eta$'s are computed on a small interval $\hat{\Omega}_i$, $i=1,2$ respectively, of size $6$.}

\begin{figure}[htbp]
\begin{center}
    \subfigure[]{\label{l13}\includegraphics[height=5.5cm]{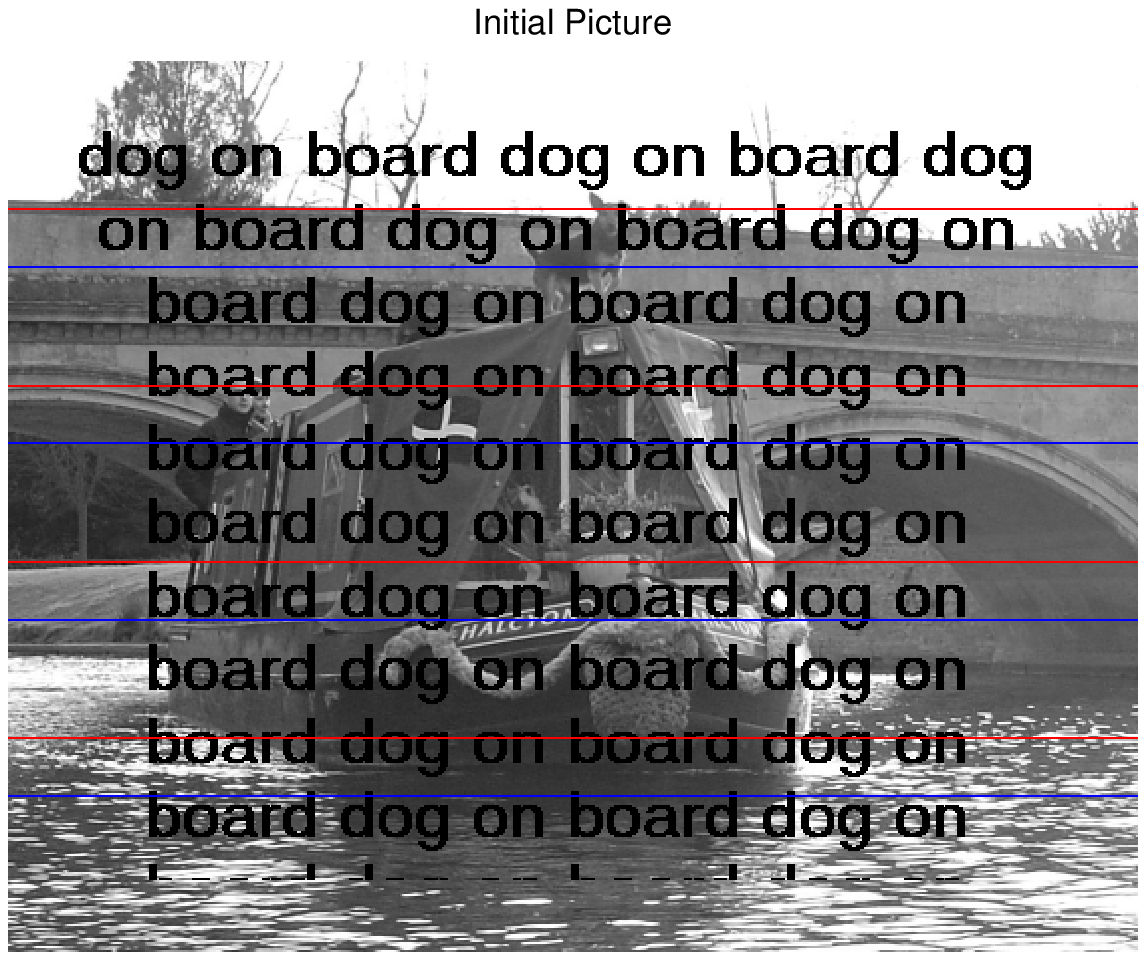}}
    \subfigure[]{\label{l23}\includegraphics[height=5.5cm]{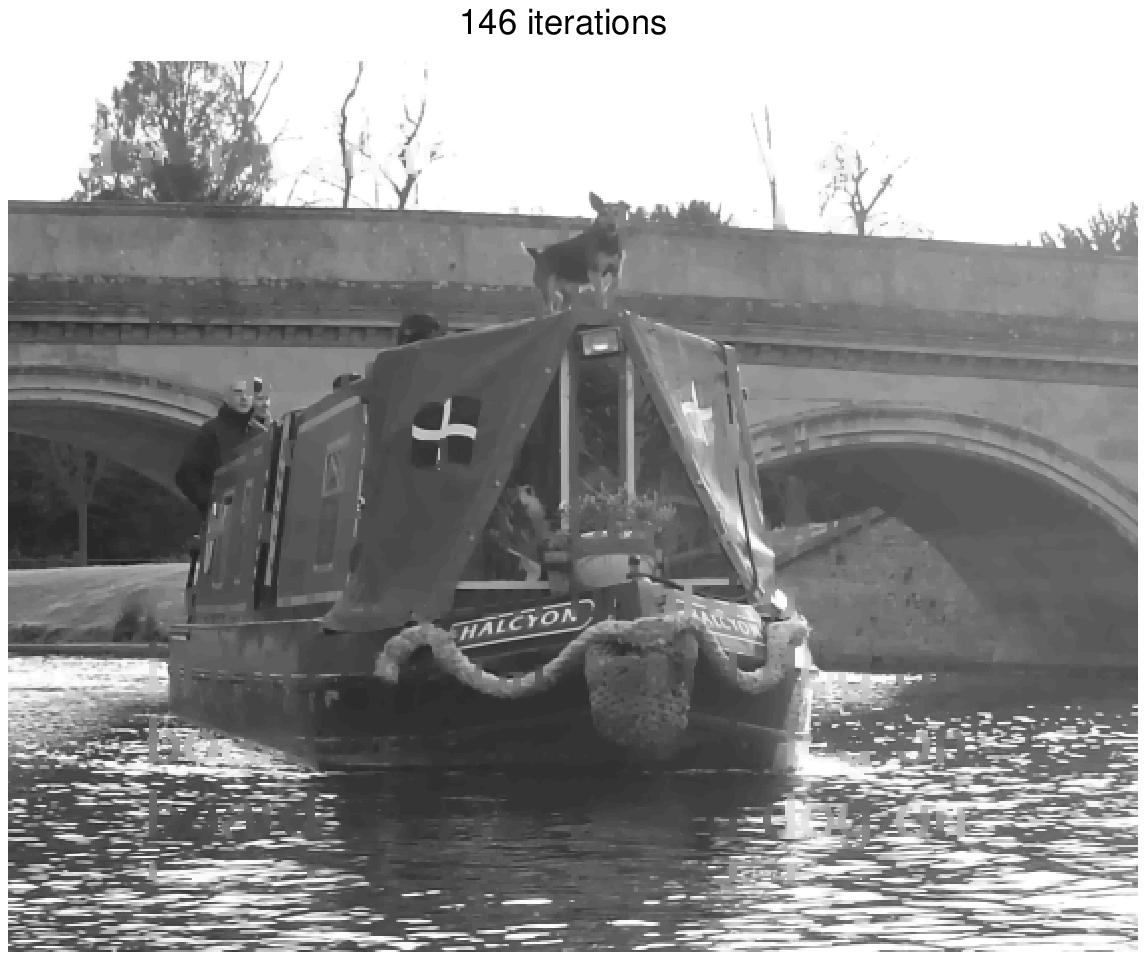}}\\
\end{center}    
\caption{\small This figure shows an application of algorithm \eqref{schw_sp:it2} for image inpainting. In this simulation the problem was split into five subproblems on overlapping subdomains.}
\label{fig:2Dinpainting}
\end{figure}

{ Figure \ref{fig:2Dinpainting} shows an example of the domain decomposition algorithm \eqref{schw_sp:it2} for total variation inpainting. As for the 1D example in Figures \ref{fig:1Dnum1}-\ref{fig:1Dnum} the operator $T$ is a multiplier, i.e., $T u =1_{\Omega\setminus D} \cdot u$, where $\Omega$ denotes the rectangular image domain and $D\subset\Omega$ the missing domain in which the original image content got lost. The regularization parameter $\alpha$ is fixed at the value $10^{-2}$. In Figure \ref{fig:2Dinpainting} the missing domain $D$ is the black writing which covers parts of the image. Here, the image domain of size $449\times 570$ pixels is split into five overlapping subdomains with an overlap size $G=28\times 570$. Further, the fixed points $\eta$'s are computed on a small stripe $\hat{\Omega}_i$, $i=1,\ldots, 5$ respectively, of size $6\times 570$ pixels.}

\begin{figure}[htbp]
\begin{center}
    \subfigure[]{\label{l11}\includegraphics[height=5cm]{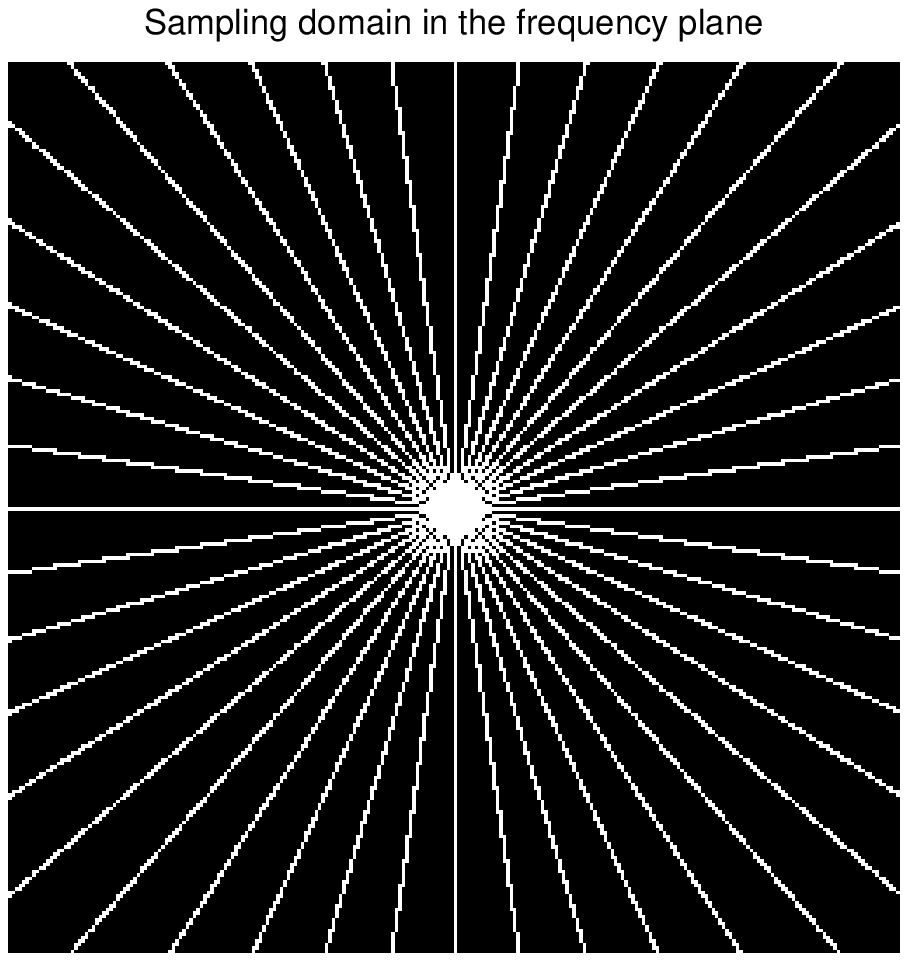}}
    \subfigure[]{\label{l21}\includegraphics[height=5cm]{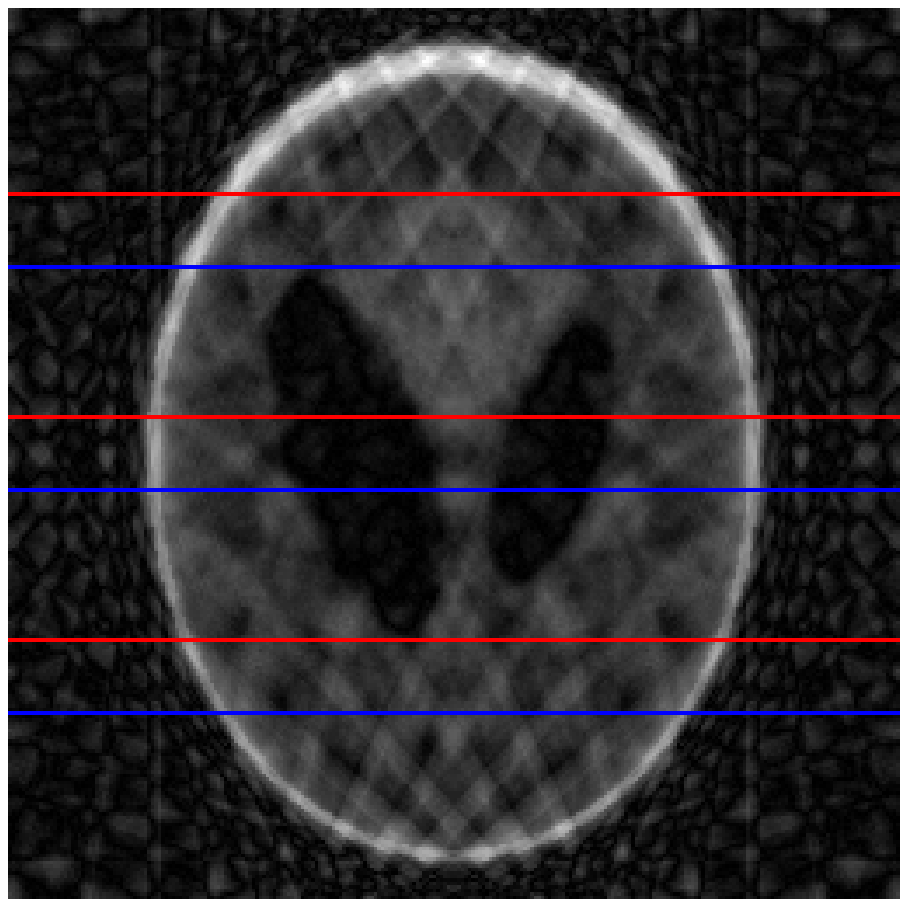}}\\
 		\subfigure[]{\label{l12}\includegraphics[height=5cm]{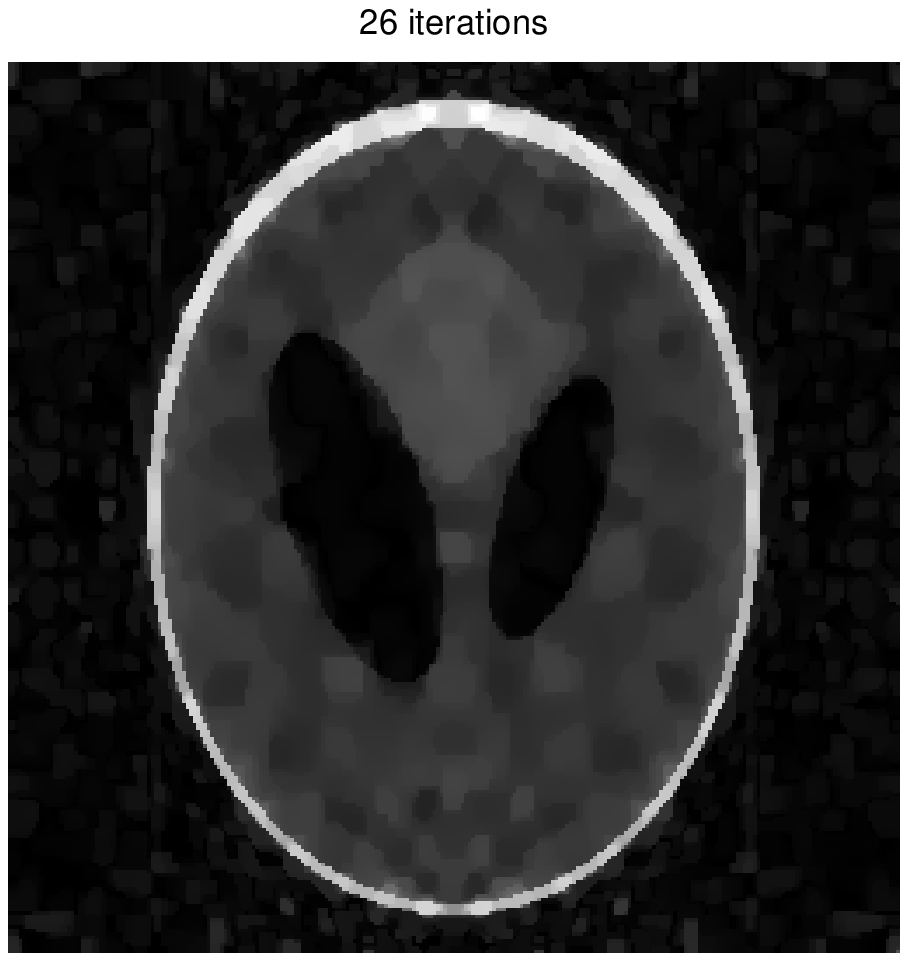}}
    \subfigure[]{\label{l22}\includegraphics[height=5cm]{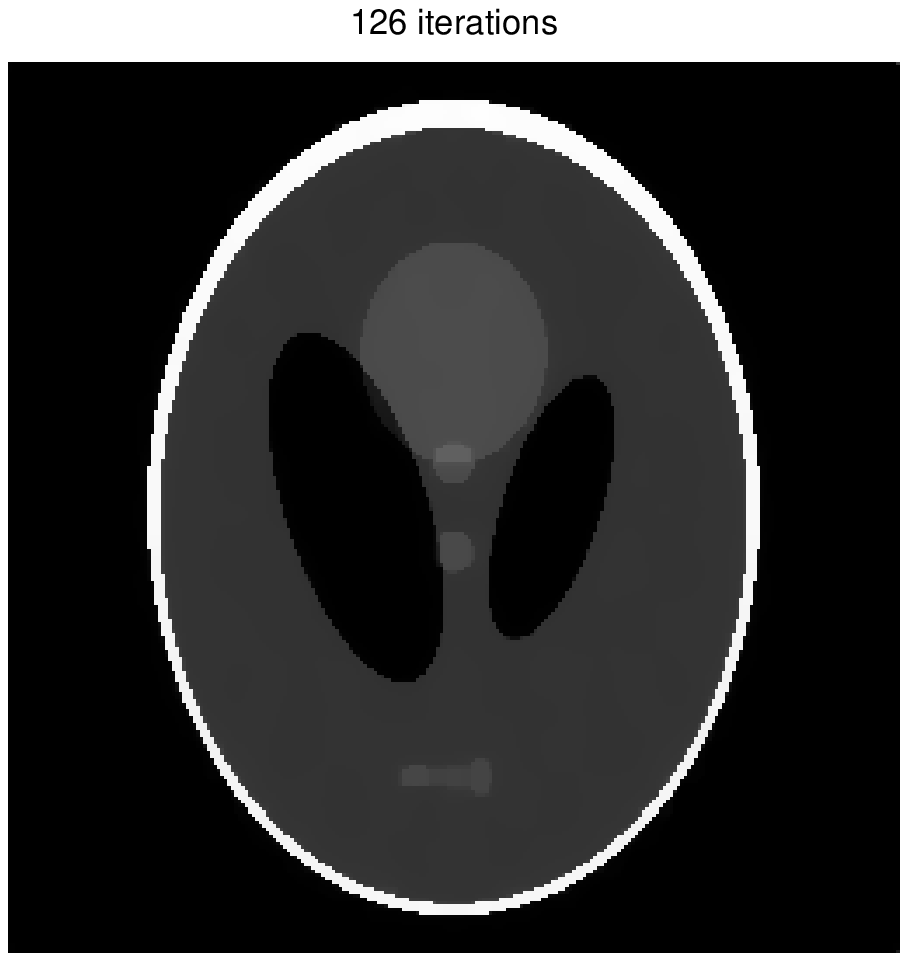}}
\end{center}    
\caption{\small We show an application of algorithm \eqref{schw_sp:it2} in a classical compressed sensing problem for recovering piecewise constant medical-type images from given partial Fourier data. In this simulation the problem was split via decomposition into four overlapping subdomains. On the top-left figure, we show the sampling data of the image in the Fourier domain. On the top-right the back-projection provided by the sampled frequency data together with the highlighted partition of the physical domain into four subdomains is shown. The bottom figures present intermediate iterations of the algorithm, i.e.,  $u^{(26)}$ and $u^{(125)}$.}
\label{fig:compressedsensing}
\end{figure}

{ Finally, in Figure \ref{fig:compressedsensing} we illustrate  the successful application of our domain decomposition algorithm \eqref{schw_sp:it2} for a compressed sensing problem. Here, we consider a medical-type image (the so-called {\it Logan-Shepp phantom}) and its reconstruction from only partial Fourier data. In this case the linear operator $T = S \circ \mathcal{F}$, where $\mathcal{F}$ denotes the $2D$ Fourier matrix and $S$ is a {\it downsampling operator} which selects only a few frequencies as output.
We minimize $\mathcal J$ with $\alpha$ set at $0.4 \times 10^{-2}$. In the application of algorithm \eqref{schw_sp:it2} the image domain of size $256\times 256$ pixels is split into four overlapping subdomains with an overlap size $G=20\times 256$. The fixed points $\eta$'s are computed in a small stripe $\hat{\Omega}_i$, $i=1,\ldots, 4$ respectively, of size $6\times 256$ pixels. 
}

 \section*{Acknowledgments}
  Massimo Fornasier and Andreas Langer acknowledge the financial support provided by the
FWF project  Y 432-N15 START-Preis {\it Sparse Approximation and Optimization in High Dimensions}. Carola-B. Sch\"onlieb acknowledges the financial support provided by the Wissenschaftskolleg (Graduiertenkolleg, Ph.D. program) of the Faculty for Mathematics at the University of Vienna (funded by the Austrian Science Fund FWF) and the FFG project no. 813610 {\it Erarbeitung neuer Algorithmen zum Image Inpainting}. Further, this publication is based on work supported by Award No. KUK-I1-007-43 , made by King Abdullah University of Science and Technology (KAUST). The results of the paper also contribute to the project WWTF Five senses-Call 2006, {\it Mathematical Methods for Image Analysis and Processing in the Visual Arts}. 
 

\appendix

\section{Proof of Proposition \ref{finitVese}}\label{AppendixB}

It is clear that $\zeta\in\partial\mathcal{J}_{\varphi}(u)$ if and only if $u=\operatorname{argmin}_{v\in\mathcal{H}}\{\mathcal{J}_{\varphi}(v) - \langle\zeta,v\rangle_{\mathcal{H}}\}$, and let us consider the following variational problem:
\renewcommand{\theequation}{$\mathcal{P}$}
\begin{equation}\label{P}
\inf_{v\in\mathcal{H}}\{\mathcal{J}_{\varphi}(v) - \langle\zeta,v\rangle_{\mathcal{H}}\} = \inf_{v\in\mathcal{H}}\{\|Tv-g\|_2^2 + 2\alpha \varphi(|\nabla v|)(\Omega) - \langle\zeta,v\rangle_{\mathcal{H}}\}
\end{equation}
We denote such an infimum by $\inf$\eqref{P}.
Now we compute \eqref{P*} the dual of \eqref{P}. Let $\mathcal{F}:\mathcal{H}\to\R$, $\mathcal{G}:\mathcal{H}\times\mathcal{H}^d \to\R$, $\mathcal{G}_1:\mathcal{H} \to\R$, $\mathcal{G}_2:\mathcal{H}^d \to\R$, such that
\begin{eqnarray*}
\mathcal{F}(v)&=&-\langle \zeta,v \rangle_{\mathcal{H}}\\
\mathcal{G}_1(w_0)&=&\|w_0-g\|_2^2\\
\mathcal{G}_2(\bar w)&=&2\alpha \varphi(|\bar w|)(\Omega)\\
\mathcal{G}(w)&=&\mathcal{G}_1(w_0) + \mathcal{G}_2(\bar w)
\end{eqnarray*}
with $w=(w_0,\bar w)\in\mathcal{H}\times\mathcal{H}^d$. Then the dual problem of \eqref{P} is given by (cf. \cite[p 60]{ET}) 
\renewcommand{\theequation}{$\mathcal{P}^*$}
\begin{equation}\label{P*}
\sup_{p^*\in\mathcal{H}\times\mathcal{H}^d}\{-\mathcal{F}^*(\Lambda^*p^*) - \mathcal{G}^*(-p^*)\}
\end{equation}
where $\Lambda:\mathcal{H}\to\mathcal{H}\times\mathcal{H}^d$ is defined by
\renewcommand{\theequation}{\arabic{section}.\arabic{equation}}
\addtocounter{equation}{-2}
$$
\Lambda v = (Tv,(\nabla v)^1,\ldots,(\nabla v)^d)
$$
and $\Lambda^*$ is its adjoint. We denote the supremum in \eqref{P*} by $\sup$\eqref{P*}. Using the definition of the conjugate function we compute $\mathcal{F}^*$ and $\mathcal{G}^*$. In particular
$$
\mathcal{F}^*(\Lambda^*p^*)=\sup_{v\in\mathcal{H}}\{ \langle \Lambda^*p^*,v \rangle_{\mathcal{H}}-\mathcal{F}(v)\}=\sup_{v\in\mathcal{H}} \langle \Lambda^*p^*+\zeta,v \rangle_{\mathcal{H}}=
\begin{cases}
0 & \Lambda^*p^*+\zeta=0\\
\infty & \text{otherwise}
\end{cases}
$$
where $p^*=(p_0^*,\bar{p}^*)$ and 
\begin{equation*}
\begin{split}
\mathcal{G}^*(p^*)&=\sup\limits_{w\in\mathcal{H}\times\mathcal{H}^d}\{\langle p^*,w \rangle_{\mathcal{H}\times\mathcal{H}^d} - \mathcal{G}(w)\}\\
&= \sup\limits_{w = (w_0,\bar w)\in\mathcal{H}\times\mathcal{H}^d}\{\langle p_0^*,w_0 \rangle_{\mathcal{H}}+\langle \bar p^*,\bar w \rangle_{\mathcal{H}^d} - \mathcal{G}_1(w_0)-\mathcal{G}_2(\bar w)\}\\
&=\sup\limits_{w_0\in\mathcal{H}}\{\langle p_0^*,w_0 \rangle_{\mathcal{H}}- \mathcal{G}_1(w_0)\}+\sup\limits_{\bar w\in\mathcal{H}^d} \{ \langle \bar p^*,\bar w \rangle_{\mathcal{H}^d} -\mathcal{G}_2(\bar w)\}\\
&=\mathcal{G}^*_1(p_0^*) + \mathcal{G}^*_2(\bar{p}^*)
\end{split}
\end{equation*}
We have that
$$
\mathcal{G}^*_1(p_0^*) = \left\langle \frac{p_0^*}{4}+g,p_0^* \right\rangle_{\mathcal{H}}
$$
and (see \cite{ET})
$$
\mathcal{G}^*_2(\bar{p}^*)=2\alpha \varphi_1^*\left(\frac{|\bar{p}^*|}{2\alpha}\right)(\Omega)
$$
if $\frac{|\bar{p}^*(x)|}{2\alpha}\in \operatorname{Dom} \varphi_1^*$, where $\varphi_1^*$ is the conjugate function of $\varphi_1$ defined by
$$
\varphi_1(s):=\varphi(|s|) \ \ s\in\R.
$$
For ease we include in Appendix \ref{AppendixC} the explicit computation of these conjugate functions. So we can write \eqref{P*} in the following way
\begin{equation}
\sup_{p^*\in\mathcal{K}}\left\{-\left\langle \frac{-p_0^*}{4}+g,-p_0^* \right\rangle_{\mathcal{H}} - 2\alpha \varphi_1^*\left(\frac{|\bar{p}^*|}{2\alpha}\right)(\Omega)\right\}
\end{equation}
where
$$
\mathcal{K}=\left\{p^*\in\mathcal{H}\times\mathcal{H}^d: \frac{|\bar{p}^*(x)|}{2\alpha}\in \operatorname{Dom} \varphi_1^* \text{ for all }x\in\Omega, \Lambda^*p^*+\zeta=0\right\}.
$$
The function $\varphi_1$ also fulfills assumption ($A_\varphi$)(ii) (i.e., there exists $c_1>0,b\geq 0$ such that $c_1 z-b\leq \varphi_1(z)\leq c_1 z+b,$ for all $z\in\R^+$). The conjugate function of $\varphi_1$ is given by $\varphi_1^*(s) = \sup_{z\in\R}\{\langle s,z \rangle-\varphi_1(z)\}$. Using the previous inequalities and that $\varphi_1$ is even (i.e., $\varphi_1(z)=\varphi_1(-z)$ for all $z\in\R$) we have
\begin{equation}\label{NN}
(\sup_{z\in\R}\{\langle s,z \rangle-c_1|z|+b\} \geq )\sup_{z\in\R}\{\langle s,z \rangle-\varphi_1(z)\} \geq \sup_{z\in\R}\{\langle s,z \rangle-c_1|z|-b\} =
\begin{cases}
-b \quad &\text{if } |s|\leq c_1\\
\infty \quad &\text{else}
\end{cases}. 
\end{equation}
In particular, one can see that $s\in\operatorname{Dom}\varphi_1^*$ if and only if $|s|\leq c_1$.

From $\Lambda^*p^*+\zeta=0$ we obtain
\begin{align*}
\langle \Lambda^*p^*,\omega\rangle_{\mathcal{H}} + \langle \zeta,\omega\rangle_{\mathcal{H}} = \langle p^*,\Lambda \omega\rangle_{\mathcal{H}^{d+1}} + \langle \zeta,\omega\rangle_{\mathcal{H}} = \langle p_0^*,T\omega \rangle_{\mathcal{H}} + \langle \bar{p}^*,\nabla \omega \rangle_{\mathcal{H}^d} + \langle \zeta,\omega \rangle_{\mathcal{H}} = 0 \ \ \text{for all } \omega\in\mathcal{H}.
\end{align*}
Then, since $\langle \bar{p}^*,\nabla \omega \rangle_{\mathcal{H}^d} = \langle-\operatorname{div}\bar{p}^*,\omega\rangle_{\mathcal{H}}$ (see Section \ref{notations}), we have
$$
T^*p_0^* - \operatorname{div}\bar{p}^* + \zeta = 0.
$$
Hence we can write $\mathcal{K}$ in the following way
$$
\mathcal{K}=\left\{p^*=(p_0^*,\bar{p}^*)\in\mathcal{H}\times\mathcal{H}^d: \frac{|\bar{p}^*(x)|}{2\alpha}\leq c_1 \text{ for all }x\in\Omega, T^*p_0^* - \operatorname{div}\bar{p}^* + \zeta = 0\right\}.
$$
We now apply the duality results from \cite[Theorem III.4.1]{ET}, since the functional in \eqref{P} is convex, continuous with respect to $\Lambda v$ in $\mathcal{H}\times\mathcal{H}^d$, and $\inf$\eqref{P} is finite. Then $\inf$\eqref{P}$=\sup$\eqref{P*}$\in\R$ and \eqref{P*} has a solution $M=(M_0,\bar{M})\in\mathcal{K}$.

Let us assume that $u$ is a solution of \eqref{P} and $M$ is a solution of \eqref{P*}. From $\inf$\eqref{P}$=\sup$\eqref{P*} we get
\begin{equation}\label{extrem-rel}
\|Tu-g\|_2^2 + 2\alpha \varphi(|\nabla u|)(\Omega) - \langle\zeta,u\rangle_{\mathcal{H}} = - \left\langle \frac{-M_0}{4}+g,-M_0 \right \rangle_{\mathcal{H}} - 2\alpha \varphi_1^*\left(\frac{|\bar{M}|}{2\alpha}\right)(\Omega)
\end{equation}
where $M=(M_0,\bar M)\in\mathcal{H}\times\mathcal{H}^d$, $\frac{|\bar{M}(x)|}{2\alpha}\leq c_1$ and $T^*M_0 - \operatorname{div}\bar{M} + \zeta = 0$, which verifies the direct implication of \eqref{cond2}. In particular 
$$
-\langle \zeta,u \rangle_{\mathcal{H}} = \langle T^*M_0,u \rangle_{\mathcal{H}} - \langle \operatorname{div}\bar{M},u \rangle_{\mathcal{H}} = \langle M_0,Tu \rangle_{\mathcal{H}} + \langle \bar{M},\nabla u \rangle_{\mathcal{H}^d},
$$
and 
\begin{equation}\label{eq3.7}
\|Tu-g\|_2^2 + \langle M_0,Tu \rangle_{\mathcal{H}} + \langle \bar{M},\nabla u \rangle_{\mathcal{H}^d} + 2\alpha \varphi(|\nabla u|)(\Omega)+ \left \langle \frac{-M_0}{4}+g,-M_0 \right\rangle_{\mathcal{H}} + 2\alpha \varphi_1^*\left(\frac{|\bar{M}|}{2\alpha}\right)(\Omega) = 0.
\end{equation}
Let us write \eqref{eq3.7} again in the following form
\begin{equation}\label{eq3.7'}
\begin{split}
\sum_{x\in\Omega} |(Tu-g)(x)|^2 + \sum_{x\in\Omega} M_0(x)(Tu)(x) + \sum_{x\in\Omega} \sum_{j=1}^{d}\bar{M}^j(x) (\nabla u)^j(x) + \sum_{x\in\Omega} 2\alpha \varphi(|(\nabla u)(x)|)& \\
+ \sum_{x\in\Omega} \left(\frac{-M_0(x)}{4}+g(x)\right)(-M_0(x)) + \sum_{x\in\Omega}2\alpha \varphi_1^*\left(\frac{|\bar{M}(x)|}{2\alpha}\right) &= 0.
\end{split}
\end{equation}
Now we have
\begin{enumerate}
\item  $2\alpha \varphi(|(\nabla u)(x)|) + \sum_{j=1}^{d}\bar{M}^j(x) (\nabla u)^j(x) + 2\alpha \varphi_1^*\left(\frac{|\bar{M}(x)|}{2\alpha}\right) \geq 2\alpha \varphi(|(\nabla u)(x)|) -\linebreak \sum_{j=1}^{d}|\bar{M}^j(x)| |(\nabla u)^j(x)| + 2\alpha \varphi_1^*\left(\frac{|\bar{M}(x)|}{2\alpha}\right)\geq 0$ by the definition of $\varphi_1^*$, since $2\alpha \varphi_1^*\left(\frac{|\bar{M}(x)|}{2\alpha}\right)=\sup_{S\in\R^d}\{ \langle \bar{M}^j(x),S \rangle_{\R^d} - 2\alpha \varphi(|S|) \}=\sup_{S\in\R^d}\{\langle |\bar{M}^j(x)|,|S| \rangle_{\R^d} - 2\alpha \varphi(|S|) \}$ .

\item $|(Tu-g)(x)|^2 + M_0(x)(Tu)(x) + \left(\frac{-M_0(x)}{4}+g(x))(-M_0(x)\right) = (((Tu)(x)-g(x)))^2 + M_0(x)((Tu)(x)-g(x)) + \left(\frac{M_0(x)}{2}\right)^2 = \left(((Tu)(x)-g(x))+\frac{M_0(x)}{2}\right)^2\geq 0$.

\end{enumerate}
Hence condition \eqref{eq3.7} reduces to
\begin{align}
&2\alpha \varphi(|(\nabla u)(x)|) + \sum_{j=1}^{d}\bar{M}^j(x) (\nabla u)^j(x) + 2\alpha \varphi_1^*\left(\frac{|\bar{M}(x)|}{2\alpha}\right)=0\quad \text{for all } x\in\Omega \label{cond1:1}\\
&-M_0(x)=2((Tu)(x)-g(x)) \quad \text{for all } x\in\Omega. \label{cond2:1}
\end{align}

Conversely, if such an $M=(M_0,\bar{M})\in\mathcal{H}\times\mathcal{H}^d$ with $\frac{|\bar{M}(x)|}{2\alpha}\leq c_1$ exists which fulfills conditions \eqref{cond1}-\eqref{cond3}, it is clear from previous considerations that equation \eqref{extrem-rel} holds. Let us denote the functional on the left side of \eqref{extrem-rel} by 
$$
P(u):=\|Tu-g\|_2^2 + 2\alpha \varphi(|\nabla u|)(\Omega) - \langle\zeta,u\rangle_{\mathcal{H}}
$$
and the functional on the right side of \eqref{extrem-rel} by 
$$
P^*(M) := - \left\langle \frac{-M_0}{4}+g,-M_0 \right \rangle_{\mathcal{H}} - 2\alpha \varphi_1^*\left(\frac{|\bar{M}|}{2\alpha}\right)(\Omega).
$$ 
We know that the functional $P$ is the functional of \eqref{P} and $P^*$ is the functional of \eqref{P*}. Hence $\inf P=\inf$\eqref{P} and $\sup P^*=\sup$\eqref{P*}. Since $P$ is convex, continuous with respect to $\Lambda u$ in $\mathcal{H}\times\mathcal{H}^d$, and $\inf$\eqref{P} is finite we know from duality results \cite[Theorem III.4.1]{ET} that $\inf$\eqref{P}$=\sup$\eqref{P*}$\in\R$. We assume that $M$ is no solution of \eqref{P*}, i.e., $P^*(M)<\sup$\eqref{P*}, and $u$ is no solution of \eqref{P}, i.e, $P(u)>\inf$\eqref{P}. Then we have that 
$$
P(u)>\inf \eqref{P}=\sup \eqref{P*}>P^*(M).
$$
Thus \eqref{extrem-rel} is valid if and only if $M$ is a solution of \eqref{P*} and $u$ is a solution of \eqref{P} which amounts to saying that $\zeta\in\partial\mathcal{J}_{\varphi}(u)$.

If additionally $\varphi$ is differentiable and $|(\nabla u)(x)|\not=0$ for $x\in\Omega$, we show that we can compute $\bar{M}(x)$ explicitly. From equation \eqref{cond1} (resp. \eqref{cond1:1}) we have
\begin{equation}\label{eq3.11}
2\alpha \varphi_1^*\left(\frac{|-\bar{M}(x)|}{2\alpha}\right) = - \langle \bar{M}(x),(\nabla u) (x)\rangle_{\R^d} - 2\alpha \varphi(|(\nabla u)(x)|). 
\end{equation}
From the definition of conjugate function we have
\begin{equation}\label{eq3.12}
\begin{split}
2\alpha \varphi_1^*\left(\frac{|-\bar{M}(x)|}{2\alpha}\right) &= 2\alpha \sup_{t\in\R}\left\{\left\langle \frac{|-\bar{M}(x)|}{2\alpha},t \right\rangle - \varphi_1(t)\right\}\\
&=2\alpha \sup_{t\geq 0}\left\{\left\langle \frac{|-\bar{M}(x)|}{2\alpha},t \right\rangle - \varphi_1(t)\right\}\\
&=2\alpha \sup_{t\geq 0}\sup_{\stackrel{S\in\R^d}{|S|=t}}\left\{\left\langle \frac{-\bar{M}(x)}{2\alpha},S \right\rangle_{\R^d} - \varphi_1(|S|)\right\}\\
&=\sup_{S\in\R^d}\left\{\left\langle -\bar{M}(x),S \right\rangle_{\R^d} - 2\alpha\varphi(|S|)(\Omega)\right\}.
\end{split}
\end{equation}
Now, if $|(\nabla u)(x)|\not=0$ for $x\in\Omega$, then it follows from \eqref{eq3.11} that the supremum is taken on in $S=|(\nabla u)(x)|$ and we have
$$
\nabla_{S}(-\langle \bar{M}(x),S \rangle_{\R^d} - 2\alpha \varphi(|S|)(\Omega))=0
$$
which implies
$$
\bar{M}^j(x)=-2\alpha \frac{\varphi'(|(\nabla u)(x)|)}{|(\nabla u)(x)|}(\nabla u)^j(x) \quad j=1,\ldots,d,
$$
and verifies \eqref{Mbar}. This finishes the proof.


\section{Computation of conjugate functions}\label{AppendixC}

Let us calculate the conjugate function of the convex function $\mathcal{G}_1(w_0)=\|w_0-g\|_2^2$. From Definition \ref{Def.conjugate} we have
$$
\mathcal{G}^*_1(p_0^*) = \sup_{w_0\in\mathcal{H}}\{\langle w_0,p_0^* \rangle_{\mathcal{H}}-\mathcal{G}_1(w_0)\}=\sup_{w_0\in\mathcal{H}}\{\langle w_0,p_0^* \rangle_{\mathcal{H}}-\langle w_0-g,w_0-g \rangle_{\mathcal{H}}\}.
$$
We set $H(w_0):=\langle w_0,p_0^* \rangle_{\mathcal{H}}-\langle w_0-g,w_0-g \rangle_{\mathcal{H}}$. To get the maximum of $H$ we calculate the G\^{a}teaux-differential at $w_0$ of $H$,
\begin{equation*}
\begin{split}
H'(w_0)&=p_0^* - 2 (w_0-g)=0
\end{split}
\end{equation*}
and we set it to zero $H'(w_0)=0$, since $H''(w_0)<0$, and we get $w_0=\frac{p_0}{2}+g$.
Thus we have that 
\begin{equation*}
\begin{split}
\sup_{w_0\in\mathcal{H}} H(w_0) = \left\langle \frac{p_0^*}{4} + g, p_0^*\right\rangle_{\mathcal{H}} = \mathcal{G}_1^*(p_0^*)
\end{split}
\end{equation*}

Now we are going to calculate the conjugate function of $\mathcal{G}_2(\bar w)=2\alpha \varphi(|\bar w|)(\Omega)$. Associated to our notations we define the space $\mathcal{H}_0^+={\R_0^+}^{N_1\times\ldots\times N_d}$. From Definition \ref{Def.conjugate} we have
\begin{equation*}
\begin{split}
\mathcal{G}_2^*(\bar{p}^*) &= \sup_{\bar{w}\in\mathcal{H}^d}\{\langle \bar{w},\bar{p}^* \rangle_{\mathcal{H}^d} - 2\alpha \varphi(|\bar{w}|)(\Omega)\}\\
&=\sup_{t\in{\mathcal{H}_0^+}} \sup_{\stackrel{\bar{w}\in\mathcal{H}^d}{|\bar{w}(x)|=t(x)}}\{\langle \bar{w},\bar{p}^* \rangle_{\mathcal{H}^d} - 2\alpha \varphi(|\bar{w}|)(\Omega)\}\\
&=\sup_{t\in{\mathcal{H}_0^+}}\{\langle t,|\bar{p}^*| \rangle_{\mathcal{H}} - 2\alpha \varphi(t)(\Omega)\}.
\end{split}
\end{equation*}
If $\varphi$ were an even function then
\begin{equation*}
\begin{split}
\sup_{t\in{\mathcal{H}_0^+}}\{\langle t,|\bar{p}^*| \rangle_{\mathcal{H}} - 2\alpha \varphi(t)(\Omega)\}&=\sup_{t\in\mathcal{H}}\{\langle t,|\bar{p}^*| \rangle_{\mathcal{H}} - 2\alpha \varphi(t)(\Omega)\}\\
&=2\alpha \sup_{t\in\mathcal{H}}\left\{\left\langle t,\frac{|\bar{p}^*|}{2\alpha} \right\rangle_{\mathcal{H}} - \varphi(t)(\Omega)\right\} \\
&= 2\alpha \varphi^*\left(\frac{|\bar{p}^*|}{2\alpha}\right)(\Omega)
\end{split}
\end{equation*}
where $\varphi^*$ is the conjugate function of $\varphi$. 

Unfortunately $\varphi$ is not even in general. To overcome this difficulty we have to choose a function which is equal to $\varphi(s)$ for $s\geq 0$ and does not change the supremum for $s<0$. For instance, one can choose $\varphi_1(s)=\varphi(|s|)$ for $s\in\R$. Then we have
\begin{equation*}
\begin{split}
\sup_{t\in{\mathcal{H}_0^+}}\{\langle t,|\bar{p}^*| \rangle_{\mathcal{H}} - 2\alpha \varphi(t)(\Omega)\}&=\sup_{t\in\mathcal{H}}\{\langle t,|\bar{p}^*| \rangle_{\mathcal{H}} - 2\alpha \varphi_1(t)(\Omega)\}\\
&=2\alpha \sup_{t\in\mathcal{H}}\left\{\left\langle t,\frac{|\bar{p}^*|}{2\alpha} \right\rangle_{\mathcal{H}} - \varphi_1(t)(\Omega)\right\} \\
&= 2\alpha \varphi_1^*\left(\frac{|\bar{p}^*|}{2\alpha}\right)(\Omega)
\end{split}
\end{equation*}
where $\varphi_1^*$ is the conjugate function of $\varphi_1$. Note that one can also choose $\varphi_1(s)=\varphi(s)$ for $s\geq 0$ and $\varphi_1(s)=\infty$ for $s<0$.

%
%

\end{document}